\documentclass[leqno,11pt]{article}  
\usepackage{amsmath}
\usepackage{amssymb}

\usepackage[german,english]{babel}
\usepackage{amsthm}

\title{On the crystalline cohomology of Deligne-Lusztig varieties}

\author{\textsc{Elmar Grosse-Kl\"onne}}
\date{}

\theoremstyle{plain} 
\newtheorem{satz}{Theorem}[section]  
\newtheorem{kor}[satz]{Corollary}  
\newtheorem{lem}[satz]{Lemma}  
 
\newcommand{\spec}{\mbox{\rm Spec}}  
\newcommand{\proj}{\mbox{\rm Proj}}
\newcommand{\quot}{\mbox{\rm Quot}}  
\newcommand{\spm}{\mbox{\rm Sp}}  
\newcommand{\spf}{\mbox{\rm Spf}}  
\newcommand{\spwf}{\mbox{\rm Spwf}}

\newcommand{\kara}{\mbox{\rm char}}  

\newcommand{\dlog}{\mbox{\rm dlog}}

\newcommand{\diag}{\mbox{\rm diag}}

\theoremstyle{remark}

\theoremstyle{definition}

\DeclareMathOperator{\Hom}{Hom}

\begin{document}
\maketitle
\footnote[0]{To look at the rigid cohomology of Deligne-Lusztig varieties was suggested to me by Adrian Iovita; useful hints I received from Marc Cabanes, Bernard Le Stum, Fabien Trihan, Sascha Orlik and Matthias Strauch: thanks to them all.}

\begin{abstract} 

Let $X\to Y^0$ be an abelian prime-to-$p$ Galois covering of smooth schemes over a perfect field $k$ of characteristic $p>0$. Let $Y$ be a smooth compactification of $Y^0$ such that $Y-Y^0$ is a normal crossings divisor on $Y$. We describe a logarithmic $F$-crystal on $Y$ whose rational crystalline cohomology is the rigid cohomology of $X$, in particular provides a natural $W[F]$-lattice inside the latter; here $W$ is the Witt vector ring of $k$. If a finite group $G$ acts compatibly on $X$, $Y^0$ and $Y$ then our construction is $G$-equivariant. As an example we apply it to Deligne-Lusztig varieties. For a finite field $k$, if ${\mathbb G}$ is a connected reductive algebraic group defined over $k$ and ${\mathbb L}$ a $k$-rational torus satisfying a certain standard condition, we obtain a meaningful equivariant $W[F]$-lattice in the cohomology ($\ell$-adic or rigid) of the corresponding Deligne-Lusztig variety and an expression of its reduction modulo $p$ in terms of equivariant Hodge cohomology groups.
 
\end{abstract}

%


\begin{center} {\bf Introduction} \end{center}

Let $k$ be a perfect field of characteristic $p>0$, let $W$ be its Witt vector ring and let $K=\quot(W)$. One of the specific interests in $p$-adic cohomology theories for $k$-varieties, as opposed to $\ell$-adic \'{e}tale cohomology $(\ell\ne p)$, lies in the hope to construct $W[F]$-lattices (i.e. $W$-lattices stable under the action of Frobenius) in the cohomology and to explicitly describe the Frobenius action on them; typically an estimate of the slopes of Frobenius on the cohomology should be given in terms of Hodge cohomology groups of proper smooth $k$-varieties. Rigid cohomology, as it stands, is a $K$-vector space valued $p$-adic cohomology theory and does not come, a priori, with natural $W[F]$-lattices {\it whose reduction modulo $p$ one could control}. If a $k$-scheme $X$ is the open complement of a smooth divisor with normal crossings on a proper smooth $k$-scheme $\overline{X}$ then its rigid cohomology $H^*_{rig}(X)$ can be computed as the rational logarithmic crystalline cohomology of $\overline{X}$ (with logarithmic poles along $\overline{X}-X$), and in this way one indeed gets a meaningful $W[F]$-lattice. However, for general smooth $X$ such $\overline{X}$ may not exist or may not be naturally at hand. Our first purpose here is to provide a computable $W[F]$-lattice in $H^*_{rig}(X)$ in the case where $X$ admits an abelian prime-to-$p$ cover $f:X\to Y^0$ to a $k$-scheme $Y^0$ which is the open complement of a smooth divisor with normal crossings on a proper smooth $k$-scheme $Y$. Namely, we construct an explicit logarithmic $F$-crystal $E$ on $Y$ whose rational logarithmic crystalline cohomology $H_{crys}^*(E,Y/W)\otimes_WK$ is identified with $H^*_{rig}(X)$; thus (the image of) $H_{crys}^*(E,Y/W)$ is a $W[F]$-lattice in $H^*_{rig}(X)$. If a finite group $G$ acts compatibly on $X$, $Y$ and $Y^0$ than it also acts on $E$ and hence on $H_{crys}^*(E,Y/W)$. 

The mere existence of a logarithmic $F$-crystal $E$ on $Y$ as above should certainly be expected for more general 
{\it tamely ramified} coverings $f:X\to Y^0$ than just abelian prime-to-$p$ coverings $f$. However, our point here is that if 
$f$ is abelian prime-to-$p$ we can {\it explicitly describe} $E$. This should be of algorithmic interest: for example in view
of the scenario described below, but we can also conceive applications to point counting algorithms for $k$-varieties 
if $k$ is finite. 

Our second purpose is to apply our construction to Deligne-Lusztig varieties. Let now $k$ be a finite field and let $\overline{k}$ denote an algebraic closure. Let ${\mathbb G}$ be a connected reductive algebraic group over $\overline{k}$, defined over $k$, and let ${\mathbb L}$ be a $k$-rational Levi subgroup of ${\mathbb G}$. A construction due to Deligne and Lusztig \cite{dellus} associates to any parabolic subgroup ${\mathbb P}$ of ${\mathbb G}$ with ${\mathbb L}$ as Levi subgroup a smooth algebraic variety $X$ over $\overline{k}$, endowed with an action of ${\mathbb G}({k})\times{\mathbb L}(k)^{opp}$, commonly referred to as the corresponding Deligne-Lusztig variety (the pioneering paper \cite{dellus} deals with the case where ${\mathbb L}$ is a $k$-rational maximal torus). The $\ell$-adic \'{e}tale cohomology with compact support ($\ell\ne p$), viewed as a virtual ${\mathbb G}({k})\times{\mathbb L}(k)^{opp}$-representation, plays a central role for the classification of all representations (in characteristic zero) of the finite group ${\mathbb G}(k)$, by now a broad and successful branch of research. The techniques developed also proved useful for the study of the representation theory of ${\mathbb G}(k)$ on vector spaces over fields of positive characteristic $\ell\ne p$ (see \cite{caen}).

Our initial observation here is that instead of $\ell$-adic \'{e}tale cohomology with compact support one may equally work with the rigid cohomology with compact support of $X$: the resulting virtual ${\mathbb G}({k})\times{\mathbb L}(k)^{opp}$-representations are the same (upon identifying the respective characteristic-zero-coefficient fields) as for $\ell$-adic \'{e}tale cohomology with compact support. If more specifically ${\mathbb L}$ is a torus satisfying a certain standard condition then it is known that $X$ can be realized as a ${\mathbb G}({k})$-equivariant Galois covering $f:X\to Y^0$ in such a way that our ${\mathbb L}(k)$ acting on $X$ becomes the Galois group for this covering, and moreover such that $Y^0$ is the open complement of a ${\mathbb G}({k})$-stable smooth divisor with normal crossings on a proper smooth $k$-scheme $Y$ with ${\mathbb G}({k})$-action. As a result our construction provides equivariant $W(\overline{k})[F]$-lattices in the cohomology of $X$ (with $W(\overline{k})$ the Witt vector ring of $\overline{k}$). In particular we obtain a geometric and explicit construction of the reduction modulo $p$ (as virtual $\overline{k}[{\mathbb G}({k})]$-modules) of the Deligne-Lusztig characters in terms of Hodge cohomology groups of equivariant vector bundles on $Y$. In fact, $X$, $Y^0$ and $Y$ are defined over $k$ and we actually get $W[F]$-lattices.

We briefly survey the content of each section. In section \ref{crys} we describe the logarithmic $F$-crystal $E$ on $Y$ mentioned at the beginning in the case where $X$ is obtained from the smooth $k$-variety $Y^0$ by adjoining a prime-to-$p$ root of a global invertible section of ${\mathcal O}_{Y^0}$ (the general case of a prime-to-$p$ abelian covering $f$ easily reduces to this one). As a logarithmic crystal, $E$ naturally decomposes into a direct sum of rank one logarithmic crystals $E(j)$, and if a finite group $G$ acts on $X$, $Y$ and $Y^0$ then the resulting action on $E$ respects this sum decomposition and is given on each $E(j)$ by suitable automorphy factors. In section \ref{relri} we analyse the overconvergent $F$-isocrystal $E^{\dagger}$ on $Y^0$ which is the push forward via $f$ of the constant $F$-isocrystal on $X$. We see that $E^{\dagger}$ is associated to $E$ by the general construction described in \cite{letr} and from this we conclude that $H_{crys}^*(E,Y/W)\otimes_WK=H^*_{rig}(X)$. Let $H^*_{rig}(X)_j$ be the direct summand of $H^*_{rig}(X)$ corresponding to the direct summand $E(j)$ of $E$. Using the general reduction-modulo-$p$ principle for crystalline cohomology from \cite{eqcr} we show that the reduction modulo $p$ of the virtual $K[G^{opp}]$-module $\sum_s(-1)^sH_{rig}^s(X)_j$ coincides with $$\sum_s(-1)^sH^s(Y,\Omega_{Y}^{\bullet}\otimes\overline{E}(j))$$where $\Omega_{Y}^{\bullet}\otimes\overline{E}(j)$ denotes the logarithmic de Rham complex of the reduction modulo $p$ of $E(j)$. An equivalent expression of this virtual $k[G^{opp}]$-module is in terms of Hodge cohomology groups for explicit equivariant vector bundles on $Y$. We also compare the rigid cohomology with the rigid cohomology with compact support. In section \ref{delufu} we prove that for general Deligne-Lusztig varieties the $\ell$-adic \'{e}tale cohomology coincides with the rigid cohomology (as virtual representations); the key argument is taken from the proof of the independence-of-$\ell$-result in \cite{dellus}. In section \ref{delugl} we look at Deligne-Lusztig varieties for ${\mathbb G}={\rm GL}_{d+1}$ (some $d\ge1$) and with ${\mathbb L}$ a maximally non-split $k$-rational torus such that ${\mathbb L}(k)={\mathbb F}_{q^{d+1}}^{\times}$. Particular attention is paid to the direct summand $E_0$ of $E$ cut out by the trivial character of ${\mathbb L}(k)={\mathbb F}_{q^{d+1}}^{\times}$: the reduction modulo $p$ of its cohomology was completely determined in \cite{holdis}. 

Let us mention that in the case ${\mathbb G}={\rm SL}_{2}$ much more substantial results have been obtained by Haastert and Jantzen \cite{haja}; they examined the crystalline cohomology of a smooth compactification of the curve $X$ itself, a method not available in higher dimensions. We hope that the present paper can be a starting point for a full generalization of the results from \cite{haja}. For general ${\mathbb G}$ and {\it tori} ${\mathbb L}$ the reductions modulo $p$ of the characters of Deligne and Lusztig have been analysed in terms of Weyl modules in \cite{ja}.\\

{\it Notations:} For a real number $r\in{\mathbb R}$ we define ${\lfloor}{r}{\rfloor}\in{\mathbb Z}$ as the integer satisfying ${\lfloor}{r}{\rfloor}\le r<{\lfloor}{r}{\rfloor}+1$. For the definition of the rigid cohomology of schemes of finite type over a perfect field $k$ with $\kara(k)>0$ we refer to \cite{berfi}. As coefficient field of the rigid cohomology of a $k$-scheme we generally use the fraction field $K$ of the ring of Witt vectors $W$ with coefficients in $k$. We denote by $\sigma$ the Frobenius endomorphism ($p$-power map) of $k$, and also the functorially induced endomorphisms of $W$ and $K$. For log crystalline cohomology we refer to \cite{kalo}.

For a polynomial ring $W[X_1,\ldots,X_n]$ in finitely many variables over $W$, Monsky and Washnitzer defined the weak formal completion $W[X_1,\ldots,X_n]^{\dagger}$ as a subalgebra of the $p$-adic completion $W[X_1,\ldots,X_n]\,\widehat{}\,$ of $W[X_1,\ldots,X_n]$. While the elements of $W[X_1,\ldots,X_n]\,\widehat{}\,$ are those power series converging on the closed unit polydisk, the elements of $W[X_1,\ldots,X_n]^{\dagger}$ are characterized as those power series satisfying a certain overconvergence condition. To a quotient algebra $A=W[X_1,\ldots,X_n]^{\dagger}/I$ (for an ideal $I\subset W[X_1,\ldots,X_n]^{\dagger}$), Meredith \cite{mer} associated a locally ringed space ${\rm Spwf}(A)$, an affine weak formal scheme. Its underlying topological space is just that of $\spec(A\otimes_Wk)$, its ring of global functions is $A$. A morphism of affine weak formal $W$-schemes $\spwf(A)\to\spwf(B)$ is smooth (resp. \'etale) if it is flat and if the induced morphism $\spec(A)\to\spec(B)$ is smooth (resp. \'etale).

 For an affine weak formal scheme $\spwf(A)$ the $K$-algebra $A\otimes_WK$ is a $K$-dagger algebra in the sense of \cite{en1dag} so that we may form the affinoid $K$-dagger space $\spm(A\otimes_WK)$. More globally, Meredith \cite{mer} defines weak formal $W$-schemes as locally ringed spaces which locally look like affine weak formal schemes. To a weak formal $W$-scheme ${\mathfrak X}$ one may associate a "generic fibre"\, ${\mathfrak X}_K$, a $K$-dagger space in the sense of \cite{en1dag}. There is a specialization map $$sp:{\mathfrak X}_K\longrightarrow{\mathfrak X}$$in the category of ringed (Grothendieck topological) spaces. The situation is completely parallel to that of formal $W$-schemes (of finite type) and $K$-rigid spaces.
  
For a flat affine weak formal $W$-scheme ${\mathfrak X}$ with special fibre the smooth affine ${k}$-scheme $X$ the rigid cohomology $H_{rig}^*(X)$ of $X$ (as defined by Berthelot, see \cite{berco}, \cite{berfi}) is the same as the de Rham cohomology of the generic fibre ${\mathfrak X}_K$ (as an affinoid $K$-dagger space) of ${\mathfrak X}$ (see \cite{en1dag} for this comparison isomorphism).

Let $\mathfrak U$ be a smooth formal $W$-scheme. We say that a closed formal subscheme ${\mathfrak D}$ of $\mathfrak U$ is a divisor with normal crossings relative to $\spf(W)$ if \'etale locally on $\mathfrak U$ the embedding of formal $W$-schemes ${\mathfrak D}\to\mathfrak U$ takes the form$$\spf(W[X_1,\ldots,X_n]\,\widehat{}\,/(X_1\cdots X_r))\longrightarrow\spf(W[X_1,\ldots,X_n]\,\widehat{}\,)$$for some $1\le r\le n$. We call ${\mathfrak D}$ a prime divisor if we can choose $r=1$. Given a finite sum ${G}=\sum_{V}b_{V}{V}$ with $b_{V}\in\mathbb{Z}$ and with prime divisors ${V}$ on $\mathfrak U$ as above, we define the invertible ${\mathcal O}_{\mathfrak U}$-module $${\mathcal O}_{\mathfrak U}(G)=\bigotimes_{V}{\mathfrak J}^{-b_{V}}_{V}$$where ${\mathfrak J}_{V}\subset{\mathcal O}_{\mathfrak U}$ is the ideal sheaf of ${V}$ in $\mathfrak U$, an invertible ${\mathcal O}_{\mathfrak U}$-module.  

\section{Equivariant $F$-Crystals}
\label{crys}

Let $k$, $W$ and $K$ be as in the introduction. Let $Y$ denote a smooth proper $k$-scheme, $Y^0\subset Y$ an open dense subscheme such that $D=Y-Y^0$ is a normal crossings divisor on $Y$ and $f:X\to Y^0$ a finite \'etale morphism of $k$-schemes. We suppose that a finite group $G$ acts compatibly (from the left) on $Y$, $Y^0$ and $X$ and that $f$ has the following form. There is a $t\in\mathbb{N}$ with $(p,t)=1$ and a unit ${{{\Pi}}}\in\Gamma(Y^0,{\mathcal O}_{Y^0})$ such that for all $g\in G$ there is an automorphy factor $\gamma_g\in\Gamma(Y,{\mathcal L}_{Y}(D))$ with ${{{\Pi}}}/g({{{\Pi}}})=\gamma^t_g$. We require that via $f$ we may identify $$X=\underline{\spec}({\mathcal O}_{Y^0}[\Xi]/(1-\Xi^t{{{\Pi}}}))$$in such a way that the lifting of the action of $G$ from $Y^0$ to $X$ is given by $g(\Xi)=\gamma_g\Xi$ for $g\in G$. (For clarification: $\underline{\spec}$ signifies {\it relative} $\spec$; we do not require that $Y^0$ or $X$ be affine.)\\

We endow $Y$ and hence its open subschemes with the log structure associated to the normal crossings divisor $D$, and we endow $\spec(k)$ with the trivial log structure; thus $Y\to \spec(k)$ is log smooth. In this section we describe a $G$-equivariant logarithmic $F$-crystal $E$ on $Y$ which in the next section will be used to compute the rigid cohomology $H^*_{rig}(X)$ of $X$ (with its $G$-action).\\

First some more notations in the characteristic-$p$-situation. We denote by ${\mathcal V}$ the set of irreducible components of $D$ (so its elements are prime divisors on $Y$). Let ${\rm div}({{{\Pi}}})$ denote the pole-zero divisor of ${{{\Pi}}}$ on $Y$: by definition this is the minimal divisor (in the usual partial ordering on the set of divisors on $Y$) with $\Pi\in{\mathcal L}_Y({\rm div}(\Pi))$. For $V\in{\mathcal V}$ we let $\mu_{V}({{{\Pi}}})\in\mathbb{Z}$ denote the multiplicity of $V$ in ${\rm div}({{{\Pi}}})$, i.e. ${\rm div}({{{\Pi}}})=\sum_{V\in{\mathcal V}}\mu_{V}({{{\Pi}}})V$. For $0\le j\le t-1$ let $$b_{V,j}=\lfloor jt^{-1}\mu_V({{{\Pi}}})\rfloor$$ and define the divisor $D(j)$ on $Y$ as $$D(j)=\sum_{V\in{\mathcal V}}b_{V,j}V.$$For $0\le j\le t-1$ we define $0\le \nu(j)\le t-1$ and $\mu(j)\in\mathbb{Z}_{\ge0}$ by requiring\begin{gather}pj=\nu(j)+\mu(j)t.\label{munude}\end{gather}Note that for any $m\in{\mathbb{Z}}$ we have$$p\lfloor jmt^{-1}\rfloor-\mu(j)m\le\lfloor\nu(j)mt^{-1}\rfloor.$$Applied to the numbers $m=\mu_{V}({{{\Pi}}})$ we therefore get\begin{gather}pD(j)-\mu(j){\rm div}({{{\Pi}}})\le D(\nu(j))\label{mudivre}\end{gather}(in the usual partial ordering on the set of divisors on $Y$).\\

Now we begin to look at liftings to characteristic $0$.\\

{\bf Definition:} A {\it local lifting datum} is a set of data $(U,{\mathfrak U},D_{\mathfrak U},{\widetilde{\Pi}},\Phi)$ as follows. $U$ is an open subscheme of $Y$ and ${\mathfrak U}$ is a lifting of $U$ to a smooth formal $W$-scheme. $D_{\mathfrak U}$ is a normal crossings divisor (relative to $\spf(W)$) on ${\mathfrak U}$ which lifts the normal crossings divisor $D\cap U$ on $U$. We have $\widetilde{\Pi}\in\Gamma({\mathfrak U},{\mathcal L}_{{\mathfrak U}}(nD_{\mathfrak U}))$ for some (unimportant) $n\in{\mathbb Z}$ and $\widetilde{\Pi}$ lifts ${{{\Pi}}}|_U$. Finally, $\Phi:{\mathfrak U}\to{\mathfrak U}$ is an endomorphism lifting the (absolute) Frobenius endomorphism of $U$, respecting $D_{\mathfrak U}$ and such that $\Phi^*:{\mathcal O}_{\mathfrak{U}}\to {\mathcal O}_{\mathfrak{U}}$ restricts to $\sigma$ on the subring $W$ of ${\mathcal O}_{\mathfrak{U}}$.

\begin{lem} $Y$ can be covered by local lifting data. 
\end{lem}

{\sc Proof:} We may cover $Y$ by open affine subschemes $U=\spec(A)$ which admit charts as follows. There exists an \'etale morphism of $k$-schemes$$\lambda:U=\spec(A)\longrightarrow\spec(k[X_1,\ldots,X_n])$$such that $D\cap U=\spec(A/\lambda^*(X_1)\cdots\lambda^*(X_r))$ for some $1\le i\le r$. We may choose a $W$-algebra $\widetilde{A}$ lifting $A$ together with an \'etale morphism of formal $\spf(W)$-schemes $$\widetilde{\lambda}:\spf(\widetilde{A})\longrightarrow\spf(W[X_1,\ldots,X_n]\,\widehat{}\,)$$lifting $\lambda$. Put ${\mathfrak U}=\spf(\widetilde{A})$. Choose elements $\lambda_1,\ldots,\lambda_r\in \widetilde{A}$ which lift $\lambda^*(X_1),\ldots,\lambda^*(X_r)$. Put $$D_{\mathfrak U}=\spf(\widetilde{A}/\widetilde{\lambda}(X_1)\cdots\widetilde{\lambda}(X_r)).$$Now ${\Pi}$ is a regular function on $Y^0$, hence $\Pi\in\Gamma(Y,{\mathcal L}_{Y}(nD))$ for some $n\in{\mathbb Z}$. Since $$\Gamma({\mathfrak U},{\mathcal L}_{{\mathfrak U}}(nD_{\mathfrak U}))\longrightarrow\Gamma(U,{\mathcal L}_{U}(n(D\cap U)))$$is surjective we may lift $\Pi|_U\in\Gamma(U,{\mathcal L}_{U}(n(D\cap U)))$ to some $\widetilde{\Pi}\in\Gamma({\mathfrak U},{\mathcal L}_{{\mathfrak U}}(nD_{\mathfrak U}))$. Finally we define the lifting $\Phi:\spf(W[X_1,\ldots,X_n]\,\widehat{}\,)\to\spf(W[X_1,\ldots,X_n]\,\widehat{}\,)$ of the Frobenius endomorphism of $\spec([X_1,\ldots,X_n])$ by $\Phi|_W=\sigma$ and $\Phi(X_i)=X_i^p$. Since $\widetilde{\lambda}$ is \'etale we get a lifting of Frobenius $\Phi:{\mathfrak U}\to{\mathfrak U}$ as desired.\hfill$\Box$\\

To define the searched for $F$-crystal $E$ on $Y$, we first define its restriction $E_U$ to $U$, for any open $U\subset Y$ which admits a local lifting datum $(U,{\mathfrak U},D_{\mathfrak U},{\widetilde{\Pi}},\Phi)$: this we can do by giving a locally free ${\mathfrak U}$-module with logarithmic connection and with a Frobenius endomorphism.\\

Let $(U,{\mathfrak U},D_{\mathfrak U},{\widetilde{\Pi}},\Phi)$ be a local lifting datum. Clearly ${\mathfrak U}\to\spf(W)$ is a log smooth lifting of $U\to \spec(k)$. We write $\Omega_{\mathfrak U}^{\bullet}$ for the logarithmic de Rham complex on ${\mathfrak U}$ with logarithmic poles along $D_{\mathfrak U}$. For $V\in{\mathcal V}$ we define a closed subscheme $V_{\mathfrak U}$ of ${\mathfrak U}$ as follows. If $V\cap U$ is empty we declare $V_{\mathfrak U}$ to be empty. Otherwise the prime divisor $V\cap U$ on $U$ lifts to a uniquely determined $W$-flat closed subscheme $V_{\mathfrak U}$ of $D_{\mathfrak U}$ (thus $V_{\mathfrak U}$ is a prime divisor on ${\mathfrak U}$ relative to $\spf(W)$). If we are given a divisor on $Y$ of the form $G=\sum_{V\in {\mathcal V}}b_VV$ with $b_V\in\mathbb Z$ then $G_{\mathfrak U}=\sum_{V}b_VV_{\mathfrak U}$ defines a lifting of the divisor $G\cap U$ on $U$. It gives rise to the line bundle ${\mathcal L}_{\mathfrak U}(G_{\mathfrak U})$ on ${\mathfrak U}(G)$ which by abuse of notation we simply denote by ${\mathcal L}_{\mathfrak U}(G)$. For example, we abusively write ${\mathcal{L}}_{\mathfrak U}({D}(j))$ instead of ${\mathcal{L}}_{{\mathfrak U}}({D}(j)_{\mathfrak U})$. 

For $0\le j\le t-1$ we define the logarithmic integrable connection$$\nabla_j:{\mathcal{L}}_{{\mathfrak U}}({D}(j))\longrightarrow {\mathcal{L}}_{{\mathfrak U}}({D}(j))\otimes_{{\mathcal O}_{{\mathfrak U}}}\Omega^{1}_{{\mathfrak U}},$$$$\quad\quad\quad\quad\quad f\quad\quad\mapsto d(f)-jt^{-1}f\dlog(\widetilde{\Pi})$$on the line bundle ${\mathcal{L}}_{{\mathfrak U}}({D}(j))$ on ${\mathfrak U}$. In the proof of \ref{crysglue} below (which deals with a more general situation) we will see that the logarithmic integrable connection $\nabla_j$ gives rise to an isomorphism $${\mathcal D}\otimes{\mathcal{L}}_{{\mathfrak U}}({D}(j))\cong {\mathcal{L}}_{{\mathfrak U}}({D}(j))\otimes{\mathcal D}$$ where ${\mathcal D}$ denotes the structure sheaf of the $p$-adically completed divided power envelope of $U$ in an exactification of the diagonal embedding $U\to{\mathfrak U}\times_{W}{\mathfrak U}$ (see the proof of \ref{crysglue} for what this means). In view of this property we conclude (see \cite{kalo}) that the module with connection $({\mathcal{L}}_{{\mathfrak U}}({D}(j)),\nabla_j)$ defines a crystal $E_U(j)$ on (the logarithmic crystalline site of) $U$ relative to $\spf(W)$. 

Now comes Frobenius; unlike the connection $\nabla_j$ it will jump between the $\mathcal{L}_{{\mathfrak U}}({D}(j))$ for various $j$. We define the map$$F:\Phi^*{\mathcal L}_{\mathfrak U}(D(j))={\mathcal L}_{\mathfrak U}(D(j))\otimes_{{\mathcal O}_{\mathfrak U},\Phi^*}{\mathcal O}_{\mathfrak U}\longrightarrow{\mathcal L}_{\mathfrak U}(pD(j)-\mu(j){\rm div}({{{\Pi}}})),$$$$\quad\quad\quad\quad\quad\quad\quad\quad\quad\quad\quad\quad f\otimes 1\quad\quad\quad\mapsto\quad\Phi^*(f)\widetilde{\Pi}^{-\mu(j)}(\frac{{\widetilde{\Pi}}^p}{\Phi^*({\widetilde{\Pi}})})^{jt^{-1}}$$($f\in {\mathcal L}_{\mathfrak U}(D(j))$). To understand this definition note that $f\in{\mathcal L}_{\mathfrak U}(D(j))$ implies $\Phi^*(f)\in{\mathcal L}_{\mathfrak U}(pD(j))$ and that $\widetilde{\Pi}^{-\mu(j)}\in{\mathcal L}_{\mathfrak U}(-\mu(j){\rm div}({{{\Pi}}}))$, and finally that ${{\widetilde{\Pi}}^p}\equiv{\Phi^*({\widetilde{\Pi}})}$ modulo $p$, hence ${{\widetilde{\Pi}}^p}/{\Phi^*({\widetilde{\Pi}})}$ is a $1$-unit and its exponentiation with $jt^{-1}$ makes sense. In view of (\ref{mudivre}) we have ${\mathcal L}_{\mathfrak U}(pD(j)-\mu(j){\rm div}({{{\Pi}}}))\subset{\mathcal L}_{\mathfrak U}(D(\nu(j)))$ so that we may and will consider $F$ as a map \begin{gather}F:\Phi^*{\mathcal L}_{\mathfrak U}(D(j))\longrightarrow {\mathcal L}_{\mathfrak U}(D(\nu(j))).\label{frogra}\end{gather}This map commutes with the connections $\nabla_j$ and $\nabla_{\nu(j)}$. This could be verified by a straightforward computation; however, in section \ref{relri} we will get this fact for free by giving a new interpretation of the collection of crystals $E_U(j)$ tensored with $K$.  

 From (\ref{frogra}) we get morphisms of crystals\begin{gather}F:\Phi^*E_U(j)\longrightarrow E_U(\nu(j))\label{frcr}\end{gather}for $0\le j\le t-1$. We define the crystal\begin{gather}E_U=\bigoplus_{0\le j\le t-1}E_U(j).\label{dece}\end{gather}and endow $E_U$ with a Frobenius structure by taking $F:\Phi^*E_U\longrightarrow E_U$ as the sum over the maps (\ref{frcr}) (each of them composed with the respective projection $\Phi^*E_U\to\Phi^*E_U(j)$ and inclusion $E_U(\nu(j))\to E_U$). We have defined an $F$-crystal $E_U$ on $U$.\\

Next we wish to glue these locally defined $F$-crystals $E_U$ to obtain a crystal $E$ on $Y$. To do this we need to check the independence of our local constructions of the chosen local lifting data. Repeating the above constructions with respect to another local lifting datum $(U',{\mathfrak U}',D'_{\mathfrak U},{\widetilde{\Pi}}',\Phi')$ we get an $F$-crystal $E'_{U'}$ on $U'$. Let $g\in G$ and let $g^*E'_{U'}$ be the pull back $F$-crystal on $g^{-1}U'$. 

\begin{lem}\label{crysglue} There is a canonical (depending on the local lifting data) isomorphism of $F$-crystals on $U\cap g^{-1}U'$,$$\beta_{g}:E_U|_{U\cap g^{-1}U'}\cong g^*E'_{U'}|_{U\cap g^{-1}U'}.$$Given a third local lifting datum over $U''\subset Y$ and another element $h\in G$, we have $\beta_{gh}=h^*(\beta_{g})\circ(\beta_{h})$ on $U\cap g^{-1}U'\cap (gh)^{-1}U''$.
\end{lem}

{\sc Proof:} Chopping of closed subschemes we may assume that $g$ induces an isomorphism $U\to U'$. We blow up ${\mathfrak U}\times_{W}{\mathfrak U}'$ along the closed formal subschemes ${V}_{{\mathfrak U}}\times {gV}_{{\mathfrak U}'}$ for $V\in{\mathcal V}$, we then remove the strict transforms of all ${\mathfrak U}\times{gV}_{{\mathfrak U}'}$ and of all ${V}_{{\mathfrak U}}\times {\mathfrak U}'$ (for all $V$) and call the result $\mathfrak{W}$. By construction, the embedding $(1,g):U\to{\mathfrak U}\times_{W}{\mathfrak U}'$ factors as canonically as \begin{gather}U\stackrel{\iota}{\longrightarrow}\mathfrak{W}\stackrel{\rho}{\longrightarrow}{\mathfrak U}\times_{W}{\mathfrak U}'.\label{exact}\end{gather}So far this is the standard construction in logarithmic crystalline cohomology: in the terminology of logarithmic geometry, $\iota$ is an exact closed embedding and $\rho$ is log \'etale so that (\ref{exact}) may be called an exactification of the closed embedding of formal log schemes $(1,g):U\to{\mathfrak U}\times_{W}{\mathfrak U}'$. In our particular situation this means that the map $u_1:\mathfrak{W}\to\mathfrak{U}$ (resp. $u_2:\mathfrak{W}\to\mathfrak{U}'$), the composite of $\rho$ with the projection to the first (resp. the second) factor of ${\mathfrak U}\times_{W}{\mathfrak U}'$, is smooth (in the classical sense), and moreover that the pull back of $D_{\mathfrak U}$ via $u_1$ is the same as the pull back of $D_{{\mathfrak U}'}$ via $u_2$: a relative normal crossings divisor on $\mathfrak{W}$ which we denote by ${D}_{\mathfrak W}$.

We view $U$ as a closed subscheme of $\mathfrak{W}$ via $\iota$. Let ${\mathcal D}$ denote the structure sheaf of the $p$-adically completed divided power envelope of $U$ in $\mathfrak{W}$; it is supported on $U$. On sufficiently small open affine pieces $T$ of $\mathfrak{W}$ the ring ${\mathcal D}(T)$ is via $u_1$ (resp. $u_2$) a relative $p$-adically completed divided power polynomial ring over ${\mathcal O}_{\mathfrak{U}}$ (resp. over ${\mathcal O}_{\mathfrak{U}'}$): this is because $u_1$ (resp. $u_2$) is smooth and admits locally a section supported on $U$. We claim that$$\alpha_j=\frac{u_1^{-1}(\widetilde{\Pi})^{jt^{-1}}}{u_2^{-1}(\widetilde{\Pi}')^{jt^{-1}}}$$is a global section and a generator (as a ${\mathcal D}$-module sheaf) of \begin{gather}{\mathcal L}_{\mathfrak{W}}(j(u_1^{-1}D(j)-u_2^{-1}D(j)))\otimes_{{\mathcal O}_{\mathfrak{W}}}{\mathcal D}.\label{dlnb}\end{gather}In particular, (\ref{dlnb}) is a free ${\mathcal D}$-module of rank one. [Here and below, to define the line bundle ${\mathcal L}_{\mathfrak{W}}(j(u_1^{-1}D(j)-u_2^{-1}D(j)))$ on ${\mathfrak{W}}$ we use similar notational conventions as we did to define line bundles on ${\mathfrak{U}}$; namely, we write the divisor $j(u_1^{-1}D(j)-u_2^{-1}D(j))$ on ${\mathfrak{W}}\otimes k$ as $\sum_{V\in{\mathcal V}}b_VV_{{\mathfrak{W}}\otimes k}$ with $b_V\in\mathbb{Z}$ and prime divisors $V_{{\mathfrak{W}}\otimes k}$ on ${\mathfrak{W}}\otimes k$. By construction, the $V_{{\mathfrak{W}}\otimes k}$ lift to $W$-flat formal subschemes $V_{{\mathfrak{W}}}$ of ${D}_{\mathfrak W}$ (which are thus prime divisors relative to $\spf(W)$), and we define the relative normal crossings divisor $j(u_1^{-1}D(j)-u_2^{-1}D(j))_{\mathfrak{W}}$ on ${\mathfrak{W}}$ as $\sum_{V\in{\mathcal V}}b_VV_{{\mathfrak{W}}}$, which we then use to define ${\mathcal L}_{\mathfrak{W}}(j(u_1^{-1}D(j)-u_2^{-1}D(j)))$.] 

We know ${{{\Pi}}}/g({{{\Pi}}})=\gamma^t_g$ for some $\gamma_g\in\Gamma(Y,{\mathcal L}_{Y}(D))$. Viewing $\gamma_g$ as an element of $\Gamma(U,{\mathcal O}_{U}(D\cap U))$ we may lift it to a section $\widetilde{\gamma}_g\in\Gamma({\mathfrak U},{\mathcal O}_{{\mathfrak U}}(D_{{\mathfrak U}}))$. Letting $g(\widetilde{\Pi})=\widetilde{\Pi}\widetilde{\gamma}_g^{-t}$ we have
\begin{gather}\alpha_j=\frac{u_1^{-1}(\widetilde{\Pi})^{jt^{-1}}}{u_1^{-1}(g(\widetilde{\Pi}))^{jt^{-1}}}\frac{u_1^{-1}(g(\widetilde{\Pi}))^{jt^{-1}}}{u_2^{-1}(\widetilde{\Pi}')^{jt^{-1}}}.\label{gammto}\end{gather}The first factor in (\ref{gammto}) is simply $u_1^{-1}(\widetilde{\gamma}_g)^j$ and this is a generator of $${\mathcal L}_{\mathfrak{W}}(j(u_1^{-1}D(j)-u_2^{-1}D(j)))={\mathcal L}_{\mathfrak{W}}(j(u_1^{-1}(D(j)-g^{-1}D(j))))$$ (note that $u_2^{-1}(V)=u_1^{-1}(g^{-1}V)$ for any divisor $V$ on $Y$). For the second factor in (\ref{gammto}) we have \begin{align}\frac{u_1^{-1}(g(\widetilde{\Pi}))^{jt^{-1}}}{u_2^{-1}(\widetilde{\Pi}')^{jt^{-1}}}&=(1+\frac{u_1^{-1}(g(\widetilde{\Pi}))-u_2^{-1}(\widetilde{\Pi}')}{u_2^{-1}(\widetilde{\Pi}')})^{jt^{-1}}\notag\\{}&=\sum_{\nu=0}^\infty\frac{jt^{-1}(jt^{-1}-1)\cdot\ldots\cdot(jt^{-1}-\nu+1)}{\nu!}(\frac{u_1^{-1}(g(\widetilde{\Pi}))-u_2^{-1}(\widetilde{\Pi}')}{u_2^{-1}(\widetilde{\Pi}')})^{\nu}.\notag\end{align}Now $\frac{u_1^{-1}(g(\widetilde{\Pi}))-u_2^{-1}(\widetilde{\Pi}')}{u_2^{-1}(\widetilde{\Pi}')}$ is a global section of the ideal of the embedding $\iota:U\to\mathfrak{W}$, hence the infinite sum is a global section and a generator of ${\mathcal D}$. (In fact the second factor in (\ref{gammto}) is an element of the structure sheaf of the completion of $\mathfrak{W}$ along $U$ since all $\frac{jt^{-1}(jt^{-1}-1)\cdot\ldots\cdot(jt^{-1}-\nu+1)}{\nu!}$ lie in ${\mathbb Z}_p$.) Our claim is established.

We have the logarithmic connection $u_1^*\nabla_j$ on $${\mathcal{L}}_{{\mathfrak W}}(u_1^{-1}D(j))=u_1^*{\mathcal{L}}_{{\mathfrak U}}({D}(j))$$ and the logarithmic connection $u_2^*\nabla_j$ on$${\mathcal{L}}_{{\mathfrak W}}(u_2^{-1}D(j))=u_2^*{\mathcal{L}}_{{\mathfrak U}'}({D}(j))$$(with logarithmic poles along ${D}_{\mathfrak W}$) and $u_1^*\nabla_j$ (resp. $u_2^*\nabla_j$) is characterized as follows. Outside ${D}_{\mathfrak W}$ we identify as usual ${\mathcal{L}}_{{\mathfrak W}}(u_1^{-1}D(j))$ (resp. ${\mathcal{L}}_{{\mathfrak W}}(u_2^{-1}D(j))$) with ${\mathcal O}_{{\mathfrak W}}$ and under this identification, $u_1^*\nabla_j$ (resp. $u_2^*\nabla_j$) acts on $f\in {\mathcal O}_{{\mathfrak W}}$ as$$u_1^*\nabla_j(f)=d(f)-jt^{-1}f{\dlog}(u_1^{-1}(\widetilde{\Pi})),$$$$u_2^*\nabla_j(f)=d(f)-jt^{-1}f{\dlog}(u_2^{-1}(\widetilde{\Pi}')).$$The connection $u_1^*\nabla_j$ (resp. the connection $u_2^*\nabla_j$) gives rise to a connection $u_1^*\nabla_j$ (resp. $u_2^*\nabla_j$) on ${\mathcal{L}}_{{\mathfrak W}}(u_1^{-1}D(j))\otimes_{{\mathcal O}_{\mathfrak W}}{\mathcal D}$ (resp. on ${\mathcal{L}}_{{\mathfrak W}}(u_2^{-1}D(j))\otimes_{{\mathcal O}_{\mathfrak W}}{\mathcal D}$). Since $\iota:U\to {\mathfrak W}$ is an exact closed embedding into a log smooth formal $W$-scheme, these data are equivalent with crystals on $U$; namely, the crystal $E_{U}(j)$ is equivalent to the one associated with $({\mathcal{L}}_{{\mathfrak W}}(u_1^{-1}D(j))\otimes_{{\mathcal O}_{\mathfrak W}}{\mathcal D},u_1^*\nabla_j)$, and the crystal $g^*E'_{U'}(j)$ is equivalent to the one associated with $({\mathcal{L}}_{{\mathfrak W}}(u_2^{-1}D(j))\otimes_{{\mathcal O}_{\mathfrak W}}{\mathcal D},u_2^*\nabla_j)$. We define the isomorphism of modules with connection $$\alpha_j:({\mathcal{L}}_{{\mathfrak W}}(u_2^{-1}D(j))\otimes_{{\mathcal O}_{\mathfrak W}}{\mathcal D},u_2^*\nabla_j)\longrightarrow({\mathcal{L}}_{{\mathfrak W}}(u_1^{-1}D(j))\otimes_{{\mathcal O}_{\mathfrak W}}{\mathcal D},u_1^*\nabla_j)$$as the one induced by multiplication with the function $\alpha_j$ defined above. It induces an isomorphism of crystals $g^*E'_{U'}(j)\cong E_{U}(j)$. Taking the sum over all $j$ gives the isomorphism of crystals $g^*E'_{U'}\cong E_{U}$; it respects the Frobenius structures, so this is in fact an isomorphism of $F$-crystals. To check the cocycle condition one has to work on the divided power envelope of (an exactification of) $(1,g,gh):U\to {\mathfrak U}\times{\mathfrak U}'\times{\mathfrak U}''$ (note that $(1,g):Y\to Y\times Y$ has the same image as $(h,gh):Y\to Y\times Y$: this is the axiom for $G$ acting (from the left) on $Y$) and there it boils down to the obvious multiplicativity of the involved automorphy factors $\alpha_j$ (one of them being the product of the other two).\hfill$\Box$\\

{\it Remark:} In fact all compatibilities (commutation and cocycle conditions between Frobenii, connections and the $G$-action, i.e. the $\alpha_j$) in the above proof, partly left out there, follow from our discussion in section \ref{relri} below where we see that the $E_U$ sit inside a $G$-equivariant overconvergent $F$-isocrystal $E^{\dagger}$ on $Y^0$ for which {\it a priori} all the compatibilities hold; from this point of view the purpose of the present section is to observe that the Frobenius, connection and $G$-action of $E^{\dagger}$ respect the integral structures $E_U$.\\

Using \ref{crysglue} with $g=h=1$ we see that we can glue our $F$-crystals $E_U$ associated with local lifting data to a global crystal $E$ on $Y$ relative to $\spf(W)$. By construction, as a crystal it is naturally decomposed as\begin{gather}E=\bigoplus_{0\le j\le t-1}E(j)\label{regen}\end{gather}such that $E_U(j)$ is the restriction of $E(j)$ to $U$. Moreover \ref{crysglue} for arbitrary $g\in G$ provides an isomorphism of $F$-crystals $g^*E\cong E$, respecting the crystal decomposition (\ref{regen}). The cocycle condition in \ref{crysglue} ensures that these isomorphisms give $E$ the structure of a $G$-equivariant $F$-crystal on $Y$. In particular, the crystalline cohomology $H_{crys}^*(Y/W,E)$ becomes a right representation of $G$.\\

Let $\Omega_{Y}^{\bullet}$ denote the logarithmic de Rham complex on $Y$. Let $\overline{E}$ resp. $\overline{E}(j)$ denote the reduction modulo $p$ of the $F$-crystal $E$, resp. the crystal $E(j)$. Thus $\overline{E}(j)$ is equivalent with the logarithmic connection $\nabla_j:{\mathcal L}_Y(D(j))\to{\mathcal L}_Y(D(j))\otimes_{{\mathcal O}_Y}\Omega_{Y}^1$. The action of $G$ on each $\overline{E}(j)$ and hence on $\overline{E}$ is then easily described follows. For $g\in G$ we need to give an isomorphism $g^*\overline{E}(j)\cong \overline{E}(j)$. It is the one corresponding to the isomorphism of line bundles with connection ${\mathcal L}_Y(g^{-1}D(j))\cong{\mathcal L}_Y(D(j))$ given by multiplication with $\gamma_g^j$.\\

{\it Remarks:} (1) If $s\in\mathbb{N}$ is such that $t$ divides $p^s-1$ then $F^s$ respects the direct sum decomposition (\ref{regen}) so that each $E(j)$ is an $F^s$-crystal.

(2) Suppose $m\in{\mathbb N}$ divides $t$ and there is a ${\Lambda}\in\Gamma(Y^0,{\mathcal O}_{Y^0})$ with ${\Lambda}^m={{{\Pi}}}$. For $0\le j\le t-1-tm^{-1}$ multiplication with local liftings $\widetilde{\Lambda}$ of ${\Lambda}$ (such that $\widetilde{\Lambda}^m=\widetilde{\Pi}$ is a lifting of ${{{\Pi}}}$ as part of a local lifting datum) induces an isomorphism of crystals\begin{gather}E(j)\cong E(j+tm^{-1}).\label{shife}\end{gather} 

We conclude this section by computing the residues of the logarithmic connections $\nabla_j$. Let $V\in{\mathcal V}$ and let $(U,{\mathfrak U},D_{\mathfrak U},{\widetilde{\Pi}},\Phi)$ be a local lifting datum such that $V_{\mathfrak U}$ is non empty. Let ${\rm Res}_V(\nabla_j)$ denote the residue along $V_{\mathfrak U}$ of the logarithmic connection $\nabla_j$ on ${\mathcal{L}}_{{\mathfrak U}}({D}(j))$.

\begin{lem}\label{resnopo}\begin{align}{\rm Res}_V(\nabla_j)&=jt^{-1}\mu_V({{{\Pi}}})-b_{V,j}\notag\\{}&=jt^{-1}\mu_V({{{\Pi}}})-\lfloor jt^{-1}\mu_V({{{\Pi}}})\rfloor.\notag\end{align}\end{lem}

{\sc Proof:} Locally on ${\mathfrak{U}}$ we find an \'{e}tale morphism of formal $W$-schemes$$\widetilde{\lambda}:{\mathfrak{U}}\longrightarrow\spf(W[X_1,\ldots,X_n]\,\widehat{}\,)$$such that $g=\widetilde{\lambda}^*(X_1)\in{\mathcal O}_{\mathfrak{U}}$ is a local equation for ${V}_{\mathfrak{U}}$ in ${\mathfrak{U}}$. Consider the composition$${\mathcal{L}}_{\mathfrak{U}}({D}(j))\otimes_{{\mathcal O}_{\mathfrak{U}}}{\mathcal O}_{{V}_{\mathfrak{U}}}\longrightarrow{\mathcal{L}}_{\mathfrak{U}}({D}(j))\otimes_{{\mathcal O}_{\mathfrak{U}}}\Omega^1_{{\mathfrak U}}\otimes_{{\mathcal O}_{\mathfrak{U}}}{\mathcal O}_{{V}_{\mathfrak{U}}}\stackrel{\beta}{\longrightarrow}{\mathcal{L}}_{\mathfrak{U}}({D}(j))\otimes_{{\mathcal O}_{\mathfrak{U}}}{\mathcal O}_{{V}_{\mathfrak{U}}}$$where the first map is the one induced by $\nabla_{j}$ and the second map $\beta$ is induced by the map ${\Omega^1_{{\mathfrak U}}}\longrightarrow{\mathcal O}_{{V}_{\mathfrak{U}}}$ which sends the class of $f\dlog(g)$ with $f\in{\mathcal O}_{{\mathfrak U}}$ to the image of $f$ in ${\mathcal O}_{{V}_{\mathfrak{U}}}$ (and vanishes on all forms which are regular along ${V}_{\mathfrak{U}}$). This composition is given by multiplication with ${\rm Res}_V(\nabla_j)$ --- by definition of ${\rm Res}_V(\nabla_j)$. To compute ${\rm Res}_V(\nabla_j)$ it is enough to evaluate the map on an arbitrary non zero element of ${\mathcal{L}}_{\mathfrak{U}}({D}(j))\otimes_{{\mathcal O}_{\mathfrak{U}}}{\mathcal O}_{{V}_{\mathfrak{U}}}$. We evaluate it on the element $g^{-b_{V,j}}\in {\mathcal{L}}_{\mathfrak{U}}({D}(j))\otimes_{{\mathcal O}_{\mathfrak{U}}}{\mathcal O}_{{V}_{\mathfrak{U}}}$. We may write $\widetilde{\Pi}=g^{-\mu_V(\Pi)}h$ with a function $h$ regular and non-zero along ${V}_{\mathfrak{U}}$. Then \begin{align}\nabla_j(g^{-b_{V,j}})&=d(g^{-b_{V,j}})-jt^{-1}g^{-b_{V,j}}\dlog(g^{-\mu_V(\Pi)}h)\notag\\{}&=g^{-b_{V,j}}\dlog(g^{-b_{V,j}})-jt^{-1}g^{-b_{V,j}}\dlog(g^{-\mu_V(\Pi)}h)\notag\\{}&=-b_{V,j}g^{-b_{V,j}}\dlog(g)+jt^{-1}\mu_V(\Pi)g^{-b_{V,j}}\dlog(g)-jt^{-1}g^{-b_{V,j}}\dlog(h)\notag\end{align}and the claim follows since $\beta$ vanishes on the term $jt^{-1}g^{-b_{V,j}}\dlog(h)$.\hfill$\Box$\\

\section{Relative rigid cohomology of tame abelian coverings}
\label{relri}

We compute the relative rigid cohomology $E^{\dagger}$ of the finite \'{e}tale morphism $f:X\to Y^0$; thus $E^{\dagger}$ will be an overconvergent $F$-isocrystal on $Y^0$, endowed with an action by $G$. Since this is a local datum we suppose that $Y^0$ is affine. Then $Y^0$ lifts to a smooth weak formal affine $W$-scheme ${\mathfrak U}=\spwf(A)$, see \cite{elk}. We choose a lifting $\widetilde{\Pi}\in A$ of ${{{\Pi}}}$ and let$$B=A[\Xi]/(1-\Xi^t\widetilde{\Pi}).$$As an $A$-module $B$ decomposes as\begin{gather}B=\bigoplus_{j=0}^{t-1}\Xi^jA.\label{torde}\end{gather}${\mathfrak X}=\spwf(B)$ is a smooth weak formal $W$-scheme lifting $X$ and the map $A\to B$ defines a lifting $\widetilde{f}:{\mathfrak X}\to{\mathfrak U}$ of $f$.\\

{\bf (I)} First let us compute $E^{\dagger}$ as an overconvergent isocrystal on $Y^0$. Let $(\Omega^{\bullet}_{A},d)$ resp. $(\Omega^{\bullet}_{B},d)$ denote the de Rham complex of $A$ resp. $B$ relative to $W$. We compute the Gauss-Manin connection \begin{gather}\nabla:B\longrightarrow B\otimes_{A}\Omega^1_{A}.\label{gama}\end{gather}In $\Omega^1_{B}$ we have $0=d(1-\Xi^t\widetilde{\Pi})=d(\Xi^t\widetilde{\Pi})$, hence $d\Xi=-t^{-1}\Xi\dlog(\widetilde{\Pi})$. It follows that$$\nabla(\Xi^j)=-jt^{-1}\Xi^j\dlog(\widetilde{\Pi})$$for any $0\le j\le t-1$. In particular we see that $\nabla$ respects the direct sum decomposition (\ref{torde}). We endow the ${\mathcal O}_{{\mathfrak U}_{\mathbb{Q}}}$-module $${\mathcal E}=(\widetilde{f}_{\mathbb{Q}})_*{\mathcal O}_{{\mathfrak X}_{\mathbb{Q}}}=\bigoplus_{j=0}^{t-1}{\mathcal O}_{{\mathfrak U}_{\mathbb{Q}}}\Xi^j$$with an (integrable overconvergent) connection $\nabla$ as follows: it respects the indicated sum decomposition, and on the $j$-summand it is given via the isomorphism of ${\mathcal O}_{{\mathfrak U}_{\mathbb{Q}}}$-modules $${\mathcal O}_{{\mathfrak U}_{\mathbb{Q}}}\cong\Xi^j{\mathcal O}_{{\mathfrak U}_{\mathbb{Q}}}, \quad 1\mapsto \Xi^j$$ (for $0\le j\le t-1$) as the connection$$\nabla_j:{\mathcal O}_{{\mathfrak U}_{\mathbb{Q}}}\to\Omega^1_{{{\mathfrak U}_{\mathbb{Q}}}},$$$$\quad\quad\quad\quad\quad\quad\quad\quad\quad\quad f\quad\mapsto d(f)-jt^{-1}f\dlog(\widetilde{\Pi}).$$This defines the overconvergent isocrystal $E^{\dagger}$ on $Y^0$ which decomposes accordingly as \begin{gather}E^{\dagger}=\bigoplus_{j=0}^{t-1}E(j)^{\dagger}.\label{torshe}\end{gather}

{\bf (II)} Next let us look at the Frobenius structure on $E^{\dagger}$. Choose an endomorphism $\Phi^*$ of the $W$-algebra $A$ which restricts to $\sigma$ on the subring $W$ and which lifts the Frobenius endomorphism of $A\otimes k$. We extend $\Phi^*$ further to an endomorphism $\Phi^*$ of $B$ lifting the $p$-power Frobenius endomorphism of $B$ by prescribing $$\Phi^*(\Xi^j)=\Xi^{\nu(j)}\widetilde{\Pi}^{-\mu(j)}(\frac{\widetilde{\Pi}^p}{\Phi^*(\widetilde{\Pi})})^{jt^{-1}}$$for $0\le j\le t-1$, with the numbers $0\le \nu(j)\le t-1$ and $\mu(j)\in\mathbb{Z}_{\ge0}$ defined through equation (\ref{munude}). It follows that the Frobenius structure $$F:\Phi^*E^{\dagger}\to E^{\dagger}$$ on $E^{\dagger}$ is given as the sum of maps$$F_j:\Phi^*E(j)^{\dagger}\to E(\nu(j))^{\dagger}$$where $F_j$ is given by the morphism of modules with connection$$({\mathcal O}_{{\mathfrak U}_{\mathbb{Q}}},\nabla_j)\otimes_{{\mathcal O}_{{\mathfrak U}_{\mathbb{Q}}},\Phi^*}{\mathcal O}_{{\mathfrak U}_{\mathbb{Q}}}\to ({\mathcal O}_{{\mathfrak U}_{\mathbb{Q}}},\nabla_{\nu(j)})$$$$\quad\quad\quad\quad\quad\quad\quad\quad\quad f\otimes 1\quad\quad\quad\quad\mapsto \Phi^*(f)\widetilde{\Pi}^{-\mu(j)}(\frac{\widetilde{\Pi}^p}{\Phi^*(\widetilde{\Pi})})^{jt^{-1}}.$$(Compare with section \ref{crys}).\\

{\bf (III)} Finally we make explicit the $G$-action on $E^{\dagger}$. Fix $g\in G$. Let $${\mathfrak Z}={\mathfrak X}\times{\mathfrak X},\quad\quad\quad\quad{\mathfrak W}={\mathfrak U}\times{\mathfrak U}$$ (products in the category of weak formal schemes), and let ${\mathfrak Z}_{K}$ (resp. ${\mathfrak W}_{K}$) denote the generic fibre (as a dagger space) of ${\mathfrak Z}$ (resp. of ${\mathfrak W}$). We view $X$ as a closed subscheme of ${\mathfrak Z}$, of ${\mathfrak U}\times{\mathfrak X}$ and of ${\mathfrak X}\times{\mathfrak U}$ via the embeddings$$X\stackrel{(1,g)}{\longrightarrow}{\mathfrak X}\times{\mathfrak X}={\mathfrak Z},$$$$X\stackrel{(f\circ 1,g)}{\longrightarrow}{\mathfrak U}\times{\mathfrak X},$$$$X\stackrel{(1,f\circ g)}{\longrightarrow}{\mathfrak X}\times{\mathfrak U},$$and we view $Y^0$ as a closed subscheme of ${\mathfrak W}$ via the embedding$$Y^0\stackrel{(1,g)}{\longrightarrow}{\mathfrak U}\times{\mathfrak U}={\mathfrak W}.$$We denote by $]X[_{{\mathfrak Z}}$ the preimage of $X$ under the specialization map $sp:{\mathfrak Z}_{K}\to {\mathfrak Z}$; thus $]X[_{{\mathfrak Z}}$ is an admissible open dagger subspace of ${\mathfrak Z}_{K}$. Similarly we define the admissible open dagger subspace $]Y^0[_{{\mathfrak W}}$ (resp. $]X[_{{\mathfrak U}\times{\mathfrak X}}$, resp. $]X[_{{\mathfrak X}\times{\mathfrak U}}$) of ${\mathfrak W}_{K}$ (resp. of $({\mathfrak U}\times{\mathfrak X})_{K}$, resp. of $({\mathfrak X}\times{\mathfrak U})_{K}$). Denote by $$u_i:]Y^0[_{{\mathfrak W}}\longrightarrow{\mathfrak W}_K\longrightarrow {\mathfrak U}_{K}$$ the projection to the $i$-th component ($i=1,2$). Since ${\mathfrak U}$ and ${\mathfrak W}$ are smooth weak formal schemes, the structure sheaf of the formal completion of ${\mathfrak W}$ along $Y^0$ is (locally) a relative formal power series ring over the structure sheaf of the $p$-adic completion of ${\mathfrak U}$. This implies (see \cite{berco}) that the morphism which $u_i$ induces on the associated rigid spaces is a fibration in relative open polydisks; by the principles of \cite{en1dag} the same is true for the morphism of dagger spaces $u_i$ itself. Moreover, with ${\mathfrak X}\to{\mathfrak U}$ also the projections ${\mathfrak X}\times{\mathfrak U}\to{\mathfrak U}\times{\mathfrak U}$ and ${\mathfrak U}\times{\mathfrak X}\to{\mathfrak U}\times{\mathfrak U}$ are finite \'{e}tale. Therefore they induce isomorphisms between the respective formal completions along $X$ and this implies --- again first for the associated rigid spaces, then by the principles of \cite{en1dag} for the dagger spaces themselves --- that the projection maps\begin{gather}]X[_{{\mathfrak U}\times{\mathfrak X}}\longleftarrow]X[_{{\mathfrak Z}}\longrightarrow ]X[_{{\mathfrak X}\times{\mathfrak U}}\label{tubis}\end{gather}are isomorphisms. Now $]X[_{{\mathfrak Z}}$ is finite \'{e}tale over $]Y^0[_{{\mathfrak W}}$, hence its relative de Rham cohomology provides a module with (integrable overconvergent) connection $({\mathcal E}',\nabla')$ on $]Y^0[_{{\mathfrak W}}$, hence an overconvergent $F$-isocrystal $(E')^{\dagger}$ on $Y^0$. On the other hand the module with connection $u_1^*({\mathcal E},\nabla)$ (resp. $u_2^*({\mathcal E},\nabla)$) on $]Y^0[_{\mathfrak W}$ corresponds to the previously defined $F$-isocrystal $E^{\dagger}$ (resp. $g^*E^{\dagger}$) on $Y^0$ by general principles of rigid cohomology (see \cite{berco}). By construction this is also the relative de Rham cohomology of $]X[_{{\mathfrak U}\times{\mathfrak X}}\to ]Y^0[_{\mathfrak W}$ (resp. of $]X[_{{\mathfrak X}\times{\mathfrak U}}\to ]Y^0[_{\mathfrak W}$). Since the maps (\ref{tubis}) are isomorphisms the canonical maps $$u_i^{-1}:u_{i}^*({\mathcal E},\nabla)\longrightarrow({\mathcal E}',\nabla')$$$(i=1,2)$ are therefore isomorphisms. Thus we may define the isomorphism of $F$-isocrystals\begin{gather}g:g^*E^{\dagger}\longrightarrow E^{\dagger}\label{giso}\end{gather}as the one corresponding to the composite $(u_1^{-1})^{-1}\circ u_2^{-1}$.

This describes the $G$-action on $E^{\dagger}$. However, we also whish to trace back the decomposition (\ref{torshe}) in this description. We keep our fixed $g\in G$ and the above notations. In view of the isomorphisms (\ref{tubis}) we can describe the ${\mathcal O}_{]Y^0[_{{\mathfrak W}}}$-module with connection $({\mathcal E}',\nabla')$ as the pull back of $({\mathcal E},\nabla)$ via $u_i$ for either $i=1$ or $i=2$, hence both $$\{u_i^{-1}(\Xi)^j\}_{0\le j\le t-1}\quad\quad (i=1,2)$$are ${\mathcal O}_{]Y^0[_{{\mathfrak W}}}$-bases. We claim that the transformation matrix between these bases is diagonal: namely, we claim that for all $j$,$$\alpha_j=\frac{u_2^{-1}(\Xi)^j}{u_1^{-1}(\Xi)^j},$$a priori an element of the fraction field of ${\mathcal O}_{{\mathfrak Z}_K}({\mathfrak Z}_K)$, is a global section of ${\mathcal O}_{]Y^0[_{{\mathfrak W}}}$. Let ${\mathcal I}$ denote the ideal in ${\mathcal O}_{{\mathfrak Z}}$ which defines the diagonal embedding $(1,1):{\mathfrak X}\to{\mathfrak Z}$. Then $\alpha_j$ is characterized by the two properties$$1-\alpha_j\in{\mathcal I},$$$$\alpha_j^{t}=\frac{u_1^{-1}(\widetilde{\Pi})^{j}}{u_2^{-1}(\widetilde{\Pi})^{j}}.$$Hence necessarily $$\alpha_j=(\frac{(u_1^{-1}(\widetilde{\Pi})}{u_2^{-1}(\widetilde{\Pi})})^{jt^{-1}}.$$Therefore it follows from our discussion in the proof of \ref{crysglue} that $\alpha_j$ is indeed a global section of ${\mathcal O}_{]Y^0[_{{\mathfrak W}}}$. By symmetry it is clear that it is even a unit. As announced it follows that the two decompositions $${\mathcal E}'=\bigoplus_{j=0}^{t-1}{\mathcal O}_{]Y^0[_{{\mathfrak W}}}u_1^{-1}(\Xi)^j=\bigoplus_{j=0}^{t-1}{\mathcal O}_{]Y^0[_{{\mathfrak W}}}u_2^{-1}(\Xi)^j$$are in fact the same and respect the connection $\nabla'$ on ${\mathcal E}'$. Moreover, the map (\ref{giso}) respects this decomposition, and its effect on the $j$-th summand can be described as follows. 

Via the isomorphism $${\mathcal O}_{{\mathfrak U}_{K}}\to{\mathcal O}_{{\mathfrak U}_{K}}\Xi^j,\quad \quad 1\mapsto \Xi^j$$the connection which $\nabla$ induces on ${\mathcal O}_{{\mathfrak U}_{K}}\Xi^j$ becomes the connection $\nabla_j$ on ${\mathcal O}_{{\mathfrak U}_{K}}$ given by$$f\mapsto d(f)-jt^{-1}f\dlog(\widetilde{\Pi}).$$Multiplication with $\alpha_j$ defines an isomorphism of ${\mathcal O}_{]Y^0[_{{\mathfrak W}}}$-modules with connection$$({\mathcal O}_{]Y^0[_{{\mathfrak W}}},u_2^*\nabla_j)\to({\mathcal O}_{]Y^0[_{{\mathfrak W}}},u_1^*\nabla_j).$$Identifying $E(j)^{\dagger}$ (resp. $g^*E(j)^{\dagger}$) with the overconvergent isocrystal associated with $({\mathcal O}_{]Y^0[_{{\mathfrak W}}},u_1^*\nabla_j)$ (resp. associated with $({\mathcal O}_{]Y^0[_{{\mathfrak W}}},u_2^*\nabla_j)$) we therefore obtain an isomorphism of overconvergent isocrystals on $Y^0$, $$g:g^*E(j)^{\dagger}\to E(j)^{\dagger}.$$

{\bf Definition:} We have $H_{rig}^*(X)=H_{rig}^*(Y^0,E^{\dagger})$. For $0\le j\le t-1$ we now define the subspace $H_{rig}^*(X)_j$ of $H_{rig}^*(X)$ as the one corresponding to the subspace $H_{rig}^*(Y^0,E^{\dagger}(j))$ of $H_{rig}^*(Y^0,E^{\dagger})$.\\

{\it Remark:} If $K^{\times}$ contains the cyclic group $T$ of all $t$-th roots of unity then we may view $f:X\to Y^0$ as the $T$-covering for which the action of $T$ on $X$ is given by $h.\Xi=h\Xi$ for all $h\in T$. Then $H_{rig}^*(X)_j$ is the subspace of $H_{rig}^*(X)$ on which $T$ acts through the character $\theta_j:T\to K^{\times}, h\mapsto h^j$. 

\begin{satz}\label{cryrig} We have a canonical $G$-equivariant and Frobenius equivariant isomorphism\begin{gather}H_{rig}^*(X)=H_{rig}^*(Y^0,E^{\dagger})\cong H_{crys}^*(Y/W,{E})\otimes_{W}K.\label{vgli}\end{gather}For $0\le j\le t-1$ it restricts to a $G$-equivariant isomorphism$$H_{rig}^*(X)_j\cong H_{crys}^*(Y/W,{E}(j))\otimes_{W}K.$$ 
\end{satz}

{\sc Proof:} By the construction described in \cite{letr} we can associate to the $F$-crystal ${E}$ on $Y$ (which is weakly non-degenerate in the terminology of \cite{letr}: it restricts to a non-degenerate $F$-crystal on the open subscheme $Y^0$ of $Y$ where the log structure is trivial) an overconvergent $F$-isocrystal on $Y^0$: but this is precisely our overconvergent $F$-isocrystal $E^{\dagger}$, as follows from the explicit descriptions of $E$ and $E^{\dagger}$ given above. These desriptions also show the coincidence of the $G$-actions. To get the isomorphism (\ref{vgli}) we proceed as in \cite{letr} 4.2, 4.4 or \cite{shiho} Corollary 2.3.9, Theorem 3.1.1, Theorem 2.4.4: the results in \cite{shiho} are formulated only for (truly) non-degenerate $F$-crystals, but as remarked in \cite{shiho} 2.4.14 they carry over to weakly non-degenerate $F$-crystals whose residues along the compactifying normal crossings divisor have no positive integers as eigenvalues; that this condition is met by our ${E}$ was checked in \ref{resnopo} (the other condition in \cite{shiho} 2.4.14 --- that the exponents of the monodromy be non-Liouville numbers --- is guaranteed by the Frobenius structure, as remarked in  \cite{letr} 4.2). That the isomorphism (\ref{vgli}) respects the $j$-parts is clear.\hfill$\Box$\\

\begin{satz}\label{brauer} For $0\le j\le t-1$ the following three virtual ${k}[G^{opp}]$-modules are the same:\\ (i) the reduction modulo $p$ of the virtual $K[G^{opp}]$-module $\sum_s(-1)^sH_{rig}^s(X)_j$\\(ii) $\sum_s(-1)^sH^s(Y,(\Omega_{Y}^{\bullet}\otimes{\mathcal L}_{Y}(D(j)),\nabla_j))$\\(iii) $\sum_{s,m}(-1)^{s+m}H^s(Y,\Omega_{Y}^{m}\otimes{\mathcal L}_{Y}(D(j)))$.
\end{satz}

{\sc Proof:} Let us describe the main result from \cite{eqcr}. Let $Y$ be a proper and smooth $k$-scheme and suppose that the finite group $G$ acts (from the right) on $Y$. Let $E$ be a locally free, finitely generated crystal of ${\mathcal O}_{Y/W}$-modules endowed with an action by $G$ (covering the action of $G$ on $Y$). For $s\in{\mathbb Z}$ let $H_{crys}^s(Y/W,E)$ denote the $s$-th crystalline cohomology group (relative to $\spf(W)$) of the crystal $E$, a finitely generated $W$-module which is zero if $s\notin[0,2\dim(Y)]$. On the other hand, the reduction modulo $p$ of the crystal $E$ is equivalent with a locally free ${\mathcal O}_Y$-module $E_k$ with connection $E_k\to E_k\otimes_{{\mathcal O}_Y}\Omega^{1}_Y$; here $\Omega^{1}_Y$ denotes the ${\mathcal O}_Y$-module of differentials of $Y/k$. Let $\Omega^{\bullet}_Y\otimes E_k$ denote the corresponding de Rham complex. The cohomology group $H^s(Y,\Omega^{\bullet}_Y\otimes E_k)$ is a finite dimensional $k$-vector space which is zero if $s\notin[0,2\dim(Y)]$. The $G$-action on $E$ provides each $H_{crys}^s(Y/W,E)$, each $H^s(Y,\Omega^{\bullet}_Y\otimes E_k)$ and each $H^s(Y,\Omega^{t}_Y\otimes E_k)$ with a $G$-action. By definition, the reduction modulo $p$ of the $K[G]$-module $H^s_{crys}(Y/W,E)\otimes_WK$ is the ${k}[G]$-module obtained by reducing modulo $p$ the $G$-stable $W$-lattice $H^s_{crys}(Y/W,E)/({\rm torsion})$ in $H^s_{crys}(Y/W,E)\otimes_WK$. Then: for any $j$, the following three virtual ${k}[G]$-modules are the same:\\ (i) the reduction modulo $p$ of the virtual $K[G]$-module $\sum_s(-1)^sH^s_{crys}(Y/W,E)\otimes_WK$\\(ii) $\sum_s(-1)^sH^s(Y,\Omega^{\bullet}_Y\otimes E_k)$\\(iii) $\sum_{s,t}(-1)^{s+t}H^s(Y,\Omega^{t}_Y\otimes E_k)$.

This result is stated for crystals in the ordinary sense, but the transposition to logarithmic crystals is immediate. Therefore, combined with Theorem \ref{cryrig}, it implies Theorem \ref{brauer}. (Of course, (ii)=(iii) is immediately clear anyway).\hfill$\Box$\\

For a collection $H^*=(H^i)_{i\in{\mathbb Z}}$ of vector spaces indexed by $\mathbb{Z}$ we write $$\chi(H^*)=\sum_i(-1)^i\dim(H^i).$$

\begin{satz}\label{maingen} (a) $\chi(H^*_{rig}(X)_{j})$ is independent of $0\le j\le t-1$.\\(b) If $0\le j\le t-1$ is such that $j\mu_V({{{\Pi}}})$ is not divisible by $t$ for all $V\in{\mathcal V}$ then$$H_{crys}^*(Y/W,{E}(j))\otimes_WK=H_{crys,c}^*(Y/W,{E}(j))\otimes_WK,$$$$H_{rig}^*(Y^0,{E}^{\dagger}_j)=H_{rig,c}^*(Y^0,{E}^{\dagger}_j).$$(c) Suppose $Y^0$ is affine and of pure dimension $d$. For all $0\le j\le t-1$ and all $m>d$ we have $H^m_{rig}(X)_{j}=0$. For $0\le j\le t-1$ as in (b) also $H^m_{rig}(X)_{j}=0$ for all $m<d$.
\end{satz}

{\sc Proof:} Clearly (c) follows from (b) by Poincar\'{e} duality. Consider the divisor $D^-(j)=D(j)-\sum_{V\in{\mathcal V}}V$ on $Y$. Just as we defined the crystal $E(j)$ departing from the divisor $D(j)$ we may now define the crystal $E^-(j)$ departing from the divisor $D^-(j)$ but using the same rule for the integrable connections $\nabla_j$. Thus if $(U,{\mathfrak U},D_{\mathfrak U},{\widetilde{\Pi}},\Phi)$ is a local lifting datum then $E(j)|_U$ is given by the connection $\nabla_j:{\mathcal L}_{\mathfrak U}(D(j))\to{\mathcal L}_{\mathfrak U}(D(j))\otimes\Omega^1_{\mathfrak U}$ while $E^-(j)|_U$ is given by the connection $\nabla_j:{\mathcal L}_{\mathfrak U}(D^-(j))\to{\mathcal L}_{\mathfrak U}(D^-(j))\otimes\Omega^1_{\mathfrak U}$. The crystalline cohomology with compact support $H_{crys,c}^*(Y/W,{E}(j))$ of $E(j)$ is just the crystalline cohomology of the subcrystal $E^-(j)$ of $E(j)$. Hence to prove (b) for crystalline cohomology we need to show that the natural map$$H_{crys}^*(Y/W,{E}^-(j))\otimes_WK\longrightarrow H_{crys}^*(Y/W,{E}(j))\otimes_WK$$is an isomorphism. To see this it is enough to show that for sufficiently small open $U\subset Y$ the map $H_{crys}^*(U/W,{E}^-(j))\otimes_WK\longrightarrow H_{crys}^*(U/W,{E}(j))\otimes_WK$ is an isomorphism. Thus we may work with a local lifting datum $(U,{\mathfrak U},D_{\mathfrak U},{\widetilde{\Pi}},\Phi)$ and need to show that the map$$H^*(U,(\Omega^{\bullet}_{{\mathfrak U}}\otimes_{{\mathcal O}_{{\mathfrak U}}}{\mathcal L}_{\mathfrak U}(D^-(j))\otimes_WK,\nabla_j))\longrightarrow H^*(U,(\Omega^{\bullet}_{{\mathfrak U}}\otimes_{{\mathcal O}_{{\mathfrak U}}}{\mathcal L}_{\mathfrak U}(D(j))\otimes_WK,\nabla_j))$$is an isomorphism. Our hypothesis on $j$ together with \ref{resnopo} implies that the residue of $({\mathcal L}_{\mathfrak U}({D}(j)),\nabla_j)$ along $V_{\mathfrak U}$ for each $V\in{\mathcal V}$ is non-zero. This means that the quotient complex$$\frac{(\Omega^{\bullet}_{{\mathfrak U}}\otimes_{{\mathcal O}_{{\mathfrak U}}}{\mathcal L}_{\mathfrak U}(D(j))\otimes_WK,\nabla_j)}{(\Omega^{\bullet}_{{\mathfrak U}}\otimes_{{\mathcal O}_{{\mathfrak U}}}{\mathcal L}_{\mathfrak U}(D^-(j))\otimes_WK,\nabla_j)}$$is acyclic, hence (b) for crystalline cohomology. But then we also get (b) for rigid cohomology from \ref{cryrig} and Poincar\'{e} duality in crystalline and rigid cohomology. Now we prove (a). By \cite{bo} we know that ${\mathbb R}\Gamma_{crys}(Y/W,{E}(j))$ is represented by a bounded complex of finitely generated free $W$-modules. Moreover we know from \cite{bo} that$${\mathbb R}\Gamma_{crys}(Y/W,{E}(j))\otimes^{\mathbb L}_{W}{k}={\mathbb R}\Gamma_{crys}(Y/{k},{E}(k)\otimes_{W}{{k}}).$$Hence\begin{align}\chi(H^*_{rig}(X)_{j})&=\chi(H^*_{crys}(Y/W,{E}(j))\otimes_WK)\notag\\{}&=\chi(H^*_{crys}(Y/{{k}},{E}(j)\otimes_{W}{{k}}))\notag\\{}&=\chi(H^*(Y,(\Omega_{Y}^{\bullet}\otimes{\mathcal L}_{Y}(D(j)),\nabla_j)))\notag\\{}&=\chi(H^*(Y,\Omega_{Y}^{\bullet}\otimes{\mathcal L}_{Y}(D(j))))\notag\end{align}where the last term is to be understood as $\sum_{s,m}(-1)^{s+t}\dim_k(H^s(Y,\Omega_{Y}^{m}\otimes{\mathcal L}_{Y}(D(j))))$. In the same sense we understand $\chi$ in the following (\ref{chil}). For any divisor $Q=\sum_{V\in{\mathcal V}}b_VV$ and any $V\in{\mathcal V}$ we claim\begin{gather}\chi(H^*(Y,\Omega_{Y}^{\bullet}\otimes{\mathcal L}_{Y}(Q)\otimes_{{\mathcal O}_Y}{\mathcal O}_V))=0.\label{chil}\end{gather}Choose a non empty open subset $U$ of $Y$ and an element $y\in{\mathcal O}_Y(U)$ which is an equation for $V\cap U$ in $U$. Then we have for each $s$ an exact sequence of sheaves on $Y$\begin{gather}0\longrightarrow{\mathcal L}_Y(Q)\otimes_{{\mathcal O}_Y}\Omega^{s-1}_V\stackrel{\wedge\dlog(y)}{\longrightarrow}{\mathcal L}_Y(Q)\otimes\Omega_Y^{s}\otimes{\mathcal O}_V\stackrel{pr}{\longrightarrow}{\mathcal L}_Y(Q)\otimes_{{\mathcal O}_Y}\Omega^{s}_V\longrightarrow0\label{fullex}\end{gather}where $pr$ denotes the natural projection map. The vanishing (\ref{chil}) follows.
Hence, if we let $Q'=Q-V$ the exact sequence$$0\longrightarrow{\mathcal L}_{Y}(Q')\longrightarrow{\mathcal L}_{Y}(Q)\longrightarrow{\mathcal L}_{Y}(Q)\otimes{\mathcal O}_V\longrightarrow0$$shows that the number $\chi(H^*(Y,\Omega_{Y}^{\bullet}\otimes{\mathcal L}_{Y}(Q)))$ is independent of the coefficients $b_V$. Statement (a) follows.\hfill$\Box$\\

{\it Remarks:} (1) Let $Y$ and $Y^0$ be as above and consider now the more general finite \'{e}tale coverings $f:X\to Y^0$ of the following form. There is a finite index set $I$ and for each $i\in I$ an element $t_i\in\mathbb{N}$ with $(p,t_i)=1$ and a unit ${{{\Pi}}}_i\in\Gamma(Y^0,{\mathcal O}_{Y^0})$ such that$$X=\underline{\spec}({\mathcal O}_{Y^0}[\Xi_i]_{i\in I}/((1-\Xi_i^{t_i}{{{\Pi}}}_i)_{i\in I})).$$Suppose again that a finite group $G$ acts compatibly on $Y$, $Y^0$ and $X$. Then all our constructions extend straightforwardly to this more general situation: the relative rigid cohomology of $f$ extends to an explicitly described $G$-equivariant logarithmic $F$-crystal $E$ on $Y^0$. As a $G$-equivariant crystal it decomposes as $$E=\bigoplus_{i\in I}\bigoplus_{0\le j_i\le t_i-1}E((j_i)_{i\in I}),$$such that each $E((j_i)_{i\in I})$ is of rank one. Theorems \ref{cryrig}, \ref{brauer} and \ref{maingen} have obvious analogs in this situation.

(2) Now assume in addition that $Y^0$ is quasiaffine: this implies that every coherent ${\mathcal O}_{Y^0}$-module is generated by its global sections (see e.g. \cite{caen} A2.10). Consider a finite \'{e}tale Galois covering $f:X\to Y^0$ with {\it abelian} Galois group $T$ of order $t$ prime to $p$ and assume that $k$ contains all ${\rm exp}(T)$-th roots of unity. We claim that $f$ has the form just described. To see this we use the action of $T$ to decompose the ${\mathcal O}_{Y^0}(Y^0)$-module ${\mathcal O}_{X}(X)$ into $t$ free direct summands of rank one. Explicitly, let $T=\prod_{i\in I}T_i$ be a decomposition into cyclic groups. For each $i\in I$ choose a character $\theta(i):T_i\to k^{\times}$ such that any other character $T_i\to k^{\times}$ is a power of $\theta(i)$. Let $\Xi_i$ denote a generator of the ${\mathcal O}_{Y^0}(Y^0)$-submodule of ${\mathcal O}_{X}(X)$ on which $T_i$ acts through $\theta(i)$ and on which all $T_s$ for $s\ne i$ act trivially. Then $\Xi_i$ is a unit in ${\mathcal O}_{Y^0}(Y^0)$; if we let ${{{\Pi}}}_i=\Xi_i^{-t_i}$ we are precisely in the situation considered above.

(3) This leads to the following natural question. Suppose $X\to Y^0$ is a finite \'{e}tale Galois covering with (possibly non-abelian) Galois group $T$ of order prime to $p$. Does the relative rigid cohomology of $f$ extend to a logarithmic $F$-crystal $E$ on $Y^0$ ? Does it so in a $G$-equivariant manner if a finite group $G$ acts compatibly on $Y$, $Y^0$ and $X$ ? We hope to return to this question in the future.

\section{The Deligne-Lusztig functor}
\label{delufu}

For the remainder of this paper we adopt the following notations. $k$ is a finite field with $q$ elements, $\overline{k}$ an algebraic closure and $W$, resp. $W(\overline{k})$, the ring of Witt vectors with coefficients in $k$, resp. in $\overline{k}$. We let $K={\quot({W})}$, $K(\overline{k})=\quot(W(\overline{k}))$ and let $\overline{K}$ denote an algebraic closure of $K$ and $K(\overline{k})$.

Let ${\mathbb G}$ be a reductive algebraic group over $\overline{k}$ with Frobenius endomorphism $F$ such that ${\mathbb G}(\overline{k})^F={\mathbb G}({k})$ (in particular ${\mathbb G}$ is $k$-rational). Let ${\mathbb L}$ be a $k$-rational Levi subgroup of ${\mathbb G}$ and let ${\mathbb P}$ denote a parabolic subgroup of ${\mathbb G}$ whose Levi subgroup is ${\mathbb L}$ but which itself is not necessarily $k$-rational. Let ${\mathbb P}={\mathbb L}{\mathbb U}$ be the Levi decomposition. Define the associated Deligne-Lusztig variety as the following subvariety of ${\mathbb G}/{\mathbb U}$:$$X=\{g.{\mathbb U}\quad;\quad g^{-1}F(g)\in {\mathbb U}.F({\mathbb U})\}.$$$X$ is a smooth (\cite{caen} Theorem 7.7) and quasiaffine (\cite{caen} Theorem 7.15) $\overline{k}$-variety and ${\mathbb G}({k})\times{\mathbb L}(k)^{opp}$ acts from the left on $X$.

Let $\ell$ be any prime number different from $p$ and fix an identification $\overline{K}\cong\overline{\mathbb Q}_{\ell}$. We define the virtual $\overline{K}[{\mathbb G}({k})\times{\mathbb L}(k)^{opp}]$-modules\begin{align}H_{et,c}^{\heartsuit}(X,\overline{\mathbb Q}_{\ell})&=\sum_i(-1)^iH_{et,c}^i(X,\overline{\mathbb Q}_{\ell}),\notag\\H_{rig,c}^{\heartsuit}(X/\overline{K})&=\sum_i(-1)^iH_{rig,c}^i(X)\otimes_{K(\overline{k})}\overline{K}\notag\end{align}and the virtual $\overline{K}[{\mathbb L}(k)\times{\mathbb G}({k})^{opp}]$-module$$H_{rig}^{\heartsuit}(X/\overline{K})=\sum_i(-1)^iH_{rig}^i(X)\otimes_{K(\overline{k})}\overline{K}.$$

The virtual $\overline{K}[{\mathbb G}({k})\times{\mathbb L}(k)^{opp}]$-module $H_{et,c}^{\heartsuit}(X,\overline{\mathbb Q}_{\ell})$ gives rise to the Deligne-Lusztig functor associated with ${\mathbb L}\subset{\mathbb G}$, see \cite{caen}.

\begin{satz}\label{ellrig} (a) $$H_{et,c}^{\heartsuit}(X,\overline{\mathbb Q}_{\ell})=H_{rig,c}^{\heartsuit}(X/\overline{K})$$ as virtual $\overline{K}[{\mathbb G}({k})\times{\mathbb L}(k)^{opp}]$-modules.\\(b) $$H_{et,c}^{\heartsuit}(X,\overline{\mathbb Q}_{\ell})=\Hom_{\overline{K}}(H_{rig}^{\heartsuit}(X/\overline{K}),{\overline{K}})$$ as virtual $\overline{K}[{\mathbb G}({k})\times{\mathbb L}(k)^{opp}]$-modules.
\end{satz}

{\sc Proof:} (a) Since the group algebra of a finite group over an algebraically closed field of characteristic $0$ is semisimple, the elements in the Grothendieck group of its modules are uniquely determined by their characters. Therefore we need to prove that for all $(g,l)\in {\mathbb G}({k})\times{\mathbb L}(k)^{opp}$ we have $$Tr((g,l)^*\,\,|\,\,H_{et,c}^{\heartsuit}(X,\overline{\mathbb Q}_{\ell}))=Tr((g,l)^*\,\,|\,\,H_{rig,c}^{\heartsuit}(X/\overline{K}))$$where $Tr$ denotes the trace. Let more generally $\sigma:X\to X$ be an automorphism of finite order. Then we claim$$Tr(\sigma^*\,\,|\,\,H_{et,c}^{\heartsuit}(X,\overline{\mathbb Q}_{\ell}))=Tr(\sigma^*\,\,|\,\,H_{rig,c}^{\heartsuit}(X/\overline{K})).$$(For the following argument compare with \cite{dellus} Theorem 3.2). For $n\ge1$ the composite $F^n\sigma$ is the Frobenius map relative to some new way of lowering the field of definition of $X$ from $\overline{k}$ to a finite subfield of $\overline{k}$ (see \cite{dimi} 3.3, 3.6). The Lefschetz fixed point formula for $\ell$-adic and for rigid cohomology (see \cite{etst} for the latter) shows\begin{gather}Tr((F^n\sigma)^*\,\,|\,\,H_{et,c}^{\heartsuit}(X,\overline{\mathbb Q}_{\ell}))=|X(\overline{k})^{F^n\sigma}|=Tr((F^n\sigma)^*\,\,|\,\,H_{rig,c}^{\heartsuit}(X/\overline{K})).\label{lefs}\end{gather}The automorphisms $F^*$ and $\sigma^*$ of the cohomology commute; hence they can be reduced to a triangular form in the same basis of $\oplus_iH_{et,c}^i(X,\overline{\mathbb Q}_{\ell})$ (resp. of $\oplus_iH_{rig,c}^i(X)\otimes_{K(\overline{k})}\overline{K}$). Hence there are $\alpha_{\lambda}$ and  $\beta_{\lambda}$ such that\begin{gather}Tr((F^n\sigma)^*\,\,|\,\,H_{et,c}^{\heartsuit}(X,\overline{\mathbb Q}_{\ell}))=\sum_{\lambda}\alpha_{\lambda}\lambda^n,\quad\quad Tr((F^n\sigma)^*\,\,|\,\,H_{rig,c}^{\heartsuit}(X/\overline{K}))=\sum_{\lambda}\beta_{\lambda}\lambda^n\label{cotr}\end{gather}for all $n\ge0$, where $\lambda$ runs through the multiplicative group of $\overline{K}\cong\overline{\mathbb Q}_{\ell}$, and where almost all $\alpha_{\lambda}$ and almost all $\beta_{\lambda}$ are zero. Comparing (\ref{lefs}) and (\ref{cotr}) for all $n\ge 1$ shows $\alpha_{\lambda}=\beta_{\lambda}$ for all $\lambda$, hence (\ref{cotr}) for $n=0$ gives our claim.

(b) This follows from (a) and Poincar\'{e}-duality in rigid cohomology since $X$ is smooth.\hfill$\Box$\\ 

Now consider the following subvariety of ${\mathbb G}/{\mathbb P}$:$$Y^0=\{g.{\mathbb P}\quad;\quad g^{-1}F(g)\in {\mathbb P}.F({\mathbb P})\}.$$Also $Y^0$ is smooth (\cite{caen} Theorem 7.7) and quasiaffine (\cite{caen} Theorem 7.15) and ${\mathbb G}(k)$ acts on it from the left. In fact $Y^0$ is the quotient of $X$ by the action of ${\mathbb L}(k)^{opp}$ so that $X\to Y^0$ is a Galois covering with group ${\mathbb L}(k)^{opp}$ (see \cite{caen} Theorem 7.8). 

The initiating paper \cite{dellus} deals with the case where ${\mathbb L}$ is a torus and ${\mathbb P}$ is a Borel subgroup. Consider the following more specific assumption. Let ${\mathbb B}_0\subset{\mathbb G}$ be a $k$-rational Borel subgroup with unipotent radical ${\mathbb U}_0$ and maximal torus ${\mathbb T}_0$. Let $S$ be the set of generating reflections of the Weyl group $W({\mathbb G},{\mathbb T}_0)=N_{{\mathbb G}}({\mathbb T}_0)/{\mathbb T}_0$ corresponding to ${\mathbb B}_0$. Let $\dot{}:W({\mathbb G},{\mathbb T}_0)\to N_{{\mathbb G}}({\mathbb T}_0)$ be the section which is multiplicative on pairs of elements whose lengths add (see \cite{caen} Theorem 7.11). Let $w\in W({\mathbb G},{\mathbb T}_0)$ be the product of some pairwise commuting elements of $S$ and let $b\in{\mathbb G}$ satisfy $b^{-1}F(b)=\dot{w}$. Then we assume that ${\mathbb U}$ is the $b$-conjugate of ${\mathbb U}_0$. (Many representation theoretic questions concerning the Deligne-Lusztig functor for general ${\mathbb L}\subset{\mathbb G}$ can be reduced to this standard assumption.) Under this assumption there is a smooth proper $\overline{k}$-scheme $Y$ with ${\mathbb G}(k)$-action and an equivariant open dense immersion $Y^0\to Y$ such that $Y-Y^0$ is a smooth divisor with normal crossings (\cite{caen} Proposition 7.13). 

Now we are precisely in the situation considered in remark (2) at the end of section \ref{relri} and our results from section \ref{relri} (extended as indicated) apply. In section \ref{delugl} below we will consider the example where ${\mathbb G}={\rm GL}_{d+1}$ and where $w$ is the longest element in the Weyl group. (In this example the Galois group $T={\mathbb L}(k)$ is cyclic, so we are in the simplified situation considered in the body of sections \ref{crys} and \ref{relri}.)

\section{Deligne-Lusztig varieties for ${\rm GL}_{d+1}$}
\label{delugl}

We fix $d\in \mathbb{N}$. Consider the affine ${k}$-scheme ${\mathbb V}$ associated with $(k^{d+1})^*=\Hom_k(k^{d+1},k)$. We write it as$${\mathbb V}=\spec({k}[\Xi_0,\ldots,\Xi_d])$$where we let $\Xi_0,\ldots,\Xi_{d}$ correspond to the canonical basis of $(k^{d+1})^*$. The right action of ${\rm GL}_{d+1}(k)={\rm GL}(k^{d+1})$ on $(k^{d+1})^*=\Hom_k(k^{d+1},k)$ defines a left action of ${\rm GL}_{d+1}(k)$ on ${\mathbb V}$.
Explicitly, ${\rm GL}_{d+1}(k)$ acts from the right on the graded ring $k[\Xi_0,\ldots,\Xi_d]$: if $f(\Xi_0,\ldots,\Xi_d)\in k[\Xi_0,\ldots,\Xi_d]$, then\begin{gather}g.f(\Xi_0,\ldots,\Xi_d)=f(\sum_{s=0}^da_{s0}\Xi_s,\ldots,\sum_{s=0}^da_{sd}\Xi_s)\quad\quad\mbox{for } g^{-1}=(a_{st})_{0\le s,t\le d}.\label{aktion}\end{gather}This induces an action from the left of ${\rm GL}_{d+1}(k)$ on ${\mathbb V}$ which passes to an action on$$Y_0={\mathbb P}((k^{d+1})^*)=\proj(k[\Xi_0,\ldots,\Xi_d])\cong {\mathbb P}_k^{d}.$$For $0\le m\le d-1$ let ${\mathcal V}_0^m$ be the set of all $k$-rational linear subvarieties $Z$ of $Y_0$ with $\dim(Z)=m$. The sequence of projective $k$-varieties$$Y=Y_{d-1}{\longrightarrow}Y_{d-2}{\longrightarrow}\ldots{\longrightarrow}Y_0$$is defined inductively by letting $Y_{m+1}\to Y_m$ be the blowing up of $Y_m$ in the strict transforms (in $Y_m$) of all $Z\in {\mathcal V}_0^m$. Let ${\mathcal V}$ denote the set of all strict transforms in $Y$ of elements of ${\mathcal V}_0^m$ for some $m$, a set of divisors on $Y$. The action of ${\rm GL}_{d+1}(k)$ on $Y_0$ naturally lifts to an action (from the left) of ${\rm GL}\sb {d+1}(k)$ on $Y$. 

On $Y_0$ and then by pull back on $Y$ we have the rational functions $z_t={\Xi_t}/{\Xi_0}$ for $0\le t\le d$. Denote by $(\Omega^{\bullet}_{Y},d)$ the de Rham complex on $Y$ with logarithmic poles along the normal crossings divisor $\sum_{V\in{\mathcal V}}V$ on $Y$. We give its open and ${\rm GL}_{d+1}(k)$-stable complement a name,$$Y^0=Y-\cup_{V\in{\mathcal V}}V.$$ 
For $0\le s\le d$ denote by ${\mathcal P_s}$ the set of subsets of $\{1,\ldots,d\}$ consisting of $s$ elements. For $1\le j\le d$ we define the rational function \begin{gather}\gamma_j=\prod_{(a_0,\ldots,a_{j-1})\in k^j}(z_j+a_{j-1}z_{j-1}+\ldots+a_1z_1+a_0)\label{gammaj}\end{gather} on $Y$, and if in addition $0\le s\le d$ we define the integer$$m_j^s=\max\{0,s-j\}q -\max\{0,s-j+1\}.$$

Recall the classification of irreducible representations of ${\rm GL}_{d+1}(k)^{opp}$ on ${k}$-vector spaces according to Carter and Lusztig. For convenience we drop the superscript $opp$ in our notation; thus, the subsequent representations are to be understood as {\it right} representations. For $1\le r\le d$ let $t_r\in{\rm GL}_{d+1}(k)$ denote the permutation matrix obtained by interchanging the $(r-1)$-st and the $r$-th row (or equivalently: column) of the identity matrix (recall that we start counting with $0$). Then $S=\{t_1,\ldots,t_d\}$ is a set of Coxeter generators for the Weyl group of ${\rm GL}_{d+1}(k)$. Let $B(k)$ resp. $U(k)$ denote the subgroup of upper triangular (resp. upper triangular unipotent) matrices.

\begin{satz} \cite{calu} (i) For an irreducible representation $\rho$ of ${\rm GL}_{d+1}(k)$ on a ${k}$-vector space, the subspace $\rho^{U(k)}$ of $U(k)$-invariants is one dimensional. If the action of $B(k)$ on $\rho^{U(k)}$ is given by the character $\chi:B(k)/U(k)\to{k}^{\times}$ and if $J=\{t\in S;\,t.\rho^{U(k)}=\rho^{U(k)}\}$, then the pair $(\chi,J)$ determines $\rho$ up to isomorphism.\\(ii) Conversely, given a character $\chi:B(k)/U(k)\to{k}^{\times}$ and a subset $J$ of $\{t\in S;\,\chi^t=\chi\}$, there exists an irreducible representation $\Theta(\chi,J)$ of ${\rm GL}_{d+1}(k)$ on a ${k}$-vector space whose associated pair (as above) is $(\chi,J)$. 
\end{satz}

In \cite{holdis} we proved the following Theorems \ref{logdif}, \ref{unimod}, \ref{bacauni}: 

\begin{satz}\label{logdif} Let $0\le s\le d$. Then $H^t(Y,\Omega_Y^s)=0$ for all $t>0$, and$$\dim_k(H^0(Y,\Omega_Y^s))=\sum_{\tau\in{\mathcal P}_s}q^{\sum_{i\in\tau}i}.$$Moreover, $H^0(Y,\Omega_Y^s)$ is generated (as a $k$-vector space) by logarithmic differential $s$-forms. 
\end{satz}

\begin{satz}\label{unimod} For $0\le s\le d$, the ${\rm GL}_{d+1}(k)$-representation on $H^0(Y,\Omega_Y^s)$ is equivalent to $\Theta(1,\{t_{s+1},\ldots,t_d\})$; it is a generalized Steinberg representation. The subspace of $U(k)$-invariants of $H^0(Y,\Omega_Y^s)$ is generated by $$\omega_s=(\prod_{j=1}^d\gamma_j^{m_j^s})dz_1\wedge\ldots\wedge dz_s.$$
\end{satz}

\begin{satz}\label{bacauni} $H_{crys}^s(Y/W)$ is torsion free for any $s$, and$$H_{crys}^s(Y/W)\otimes_{W}k=H^s(Y,\Omega_Y^{\bullet})=H^0(Y,\Omega_Y^{s}).$$
\end{satz}

Let ${\mathbb L}$ be a maximally non split torus in ${\rm GL}_{d+1}$. Thus ${\mathbb L}$ is in relative position $w$ to the standard torus ${\mathbb T}_0$ of diagonal matrices where $w$ denotes the permutation matrix with entry $1$ at position $(i,d-i)$ for each $0\le i\le d$. It follows that $T={\mathbb L}(k)$ is (abstractly) isomorphic with the set of fixed points of $wF$ acting on ${\mathbb T}_0(\overline{k})$ if $F$ denotes the standard Frobenius which raises each matrix entry to its $q$-th power. But this is the subgroup of ${\rm GL}_{d+1}(\overline{k})$ consisting of the diagonal matrices $\diag(t,t^q,\ldots,t^{d^d})$ for $t\in {\mathbb F}^{\times}_{q^{d+1}}$. Thus we may henceforth identify $T$ with the multiplicative group of the field ${\mathbb F}_{q^{d+1}}$ with $q^{d+1}$ elements. If we let $$\delta=\det((\Xi_i^{q^j})_{0\le i,j\le d})\in\mathbb{Z}[\Xi_0,\ldots,\Xi_d]$$then it is straightforwardly checked that in $k[\Xi_0,\ldots,\Xi_d]$ we have\begin{gather}\prod_{a\in k^{d+1}-\{0\}}\sum_{i=0}^da_i\Xi_i=(-1)^{d+1}\delta^{q-1}.\label{lemdet}\end{gather}

The Deligne-Lusztig variety as defined in section \ref{delufu} is the $\overline{k}$-variety $$X_{\overline{k}}=\spec(B),$$$$B=\overline{k}[\Xi_0,\ldots,\Xi_d]/(\delta^{q-1}-(-1)^d).$$In fact, only for even $d$ this is the explicit formula \cite{dellus} (2.2.2); however, since (for all $d$) we have an isomorphism of $\overline{k}$-algebras$$B\cong\overline{k}[\Xi_0,\ldots,\Xi_d]/(\delta^{q-1}-(-1)^{d+1})$$(if $\xi$ is a $(q^{d+1}-1)$-st root of $-1$ in $\overline{k}$, send $\Xi_i$ to $\xi\Xi_i$, for $i=0,\ldots,d$) our formula is equivalent with \cite{dellus} (2.2.2) also for odd $d$. This argument also shows that for {\it any} $d$ our $X_{\overline{k}}$ decomposes into $q-1$ connected components, according to the decomposition $\delta^{q-1}-1=\prod_{a\in k^{\times}}(\delta+a)$. The action of ${\rm GL}_{d+1}(k)$ on ${\mathbb V}$ passes to an action of ${\rm GL}_{d+1}(k)$ on the closed subscheme $X_{\overline{k}}$ of ${\mathbb V}\otimes_k\overline{k}$. Let $z_1,\ldots,z_d$ denote free variables, let $z_0=1$ and set$${{{\Pi}}}=-\prod_{a\in k^{d+1}-\{0\}}\sum_{i=0}^da_iz_i,$$$$A=\overline{k}[z_1,\ldots,z_d][{{{\Pi}}}^{-1}],$$$$Y^0_{\overline{k}}=\spec(A).$$$Y^0_{\overline{k}}$ is the complement in ${\mathbb P}_{\overline{k}}^d$ of all $k$-rational linear hyperplanes. Sending $z_i$ to $\Xi_i/\Xi_0$ we have an isomorphism (use (\ref{lemdet}))$$B=A[\Xi_0]/(1-\Xi_0^{q^{d+1}-1}{{{\Pi}}}).$$In this way we may view $A$ as a subring of $B$; as such it is stable under the action by ${\rm GL}_{d+1}(k)$. In fact, $X_{\overline{k}}\to Y^0_{\overline{k}}$ is a ${\rm GL}_{d+1}(k)$-equivariant finite \'{e}tale Galois covering with group $T$ if we let $T$ act on $B$ as follows: on the subring $A$ of $B$ it acts trivially, and on the class of $\Xi_0$ it acts by multiplication, $h.\Xi_0=h\Xi_0$ for $h\in T$. For $0\le j\le q^{d+1}-2$ we introduce the characters $$\theta_j:T\to K(\overline{k})^{\times},\quad h\mapsto h^j$$ (via the Teichm\"uller lifting).

Using the same defining equations we see that $X_{\overline{k}}\to Y^0_{\overline{k}}$ is obtained by base change from a morphism of $k$-schemes $X\to Y^0$. We have the proper smooth $k$-scheme $Y$ with ${\rm GL}_{d+1}(k)$-action and an equivariant open immersion $Y^0\to Y$ such that $D=Y-Y^0$ is a normal crossings divisor on $Y$. Thus all the results from sections \ref{crys} and \ref{relri} are available for $X_{\overline{k}}\to Y^0_{\overline{k}}$ and for $X\to Y^0$ (however the $T$-action does not descend from $X_{\overline{k}}$ to $X$). In fact since it is known that rigid cohomology commutes with the base change $k\to\overline{k}$ it follows that the $K(\overline{k})$-vector spaces $H_{rig}^*(X_{\overline{k}})$ and their $\theta_j$-eigenspaces $H_{rig}^*(X_{\overline{k}})_j$ are obtained by base extension $K\to K(\overline{k})$ from the $K$-vector spaces $H_{rig}^*(X)$ and their subspaces $H_{rig}^*(X)_j$. We point out that the automorphy factors $\gamma_g\in\Gamma(Y,{\mathcal L}_Y(D))$ for $g\in{\rm GL}_{d+1}(k)$ considered in section \ref{crys} are in this case given by$$\gamma_g=\frac{g(\Xi_0)}{\Xi_0}.$$The condition from \ref{maingen} (b) that $0\le j\le t-1$ be such that $j\mu_V({{{\Pi}}})$ is not divisible by $t$ for all $V\in{\mathcal V}$ becomes the condition that $j$ be not divisible by $\sum_{i=1}^dq^i=(q-1)^{-1}(q^{d+1}-1)$.  The $\theta_j$ for such $j$ are called {\it non singular} or {\it in general position}; equivalently, $\theta_j$ does not factor through a norm map $T={\mathbb F}_{q^{d+1}}^{\times}\to {\mathbb F}_{q^{s}}^{\times}$ for some $s<d+1$.

Since ${{{\Pi}}}$ is the $(q-1)$-st power of an element in $A$ the isomorphism (\ref{shife}) yields isomorphisms $$H_{rig}^*(X_{\overline{k}})_j\cong H_{rig}^*(X_{\overline{k}})_{j+\sum_{i=0}^dq^i}.$$In particular, if $\theta_j$ is singular (i.e. not in general position) then\begin{gather}H_{rig}^*(X_{\overline{k}})_j\cong H_{rig}^*(X_{\overline{k}})_{0}=H_{rig}^*(Y^0_{\overline{k}}).\label{spepo}\end{gather}

\begin{kor}\label{dimfo} The $\theta_j$-eigenspaces for non singular $\theta_j$ in the (compactly supported) rigid and $\ell$-adic \'{e}tale cohomology groups (any $\ell\ne p$) of $X_{\overline{k}}$ are non zero only in degree $d$, and in degree $d$ they coincide, as ${\rm GL}_{d+1}(k)$-representations on characteristic zero vector spaces of dimension\begin{align}\dim_{K(\overline{k})}H^d_{rig}(X_{\overline{k}})_{j}&=\sum_{s\ge0}(-1)^{s+d}\sum_{\tau\in{\mathcal P}_s}q^{\sum_{i\in\tau}i}\notag\\{}&=(q-1)(q^2-1)\cdots(q^d-1).\notag\end{align}Similarly, the $\theta_0$-eigenspaces in (compactly supported) rigid and $\ell$-adic \'{e}tale cohomology coincide as ${\rm GL}_{d+1}(k)$-representations, individually in each degree.
\end{kor}

{\sc Proof:} The coincidence of rigid and $\ell$-adic \'{e}tale cohomology as virtual ${\rm GL}_{d+1}(k)\times T$-representations follows from \ref{ellrig}. We saw in \ref{maingen} that the $\theta_j$-eigenspaces in (compactly supported) rigid cohomology live only in degree $d$, and in \cite{dellus} the same is shown for (compactly supported) $\ell$-adic cohomology. Since everything is semisimple we get coincidence. By the proof of \ref{maingen} (a) the vector space dimension is given for {\it any} $j$ by $\chi(H^*(Y,\Omega_Y^{\bullet}))$; but this is computed in \ref{logdif}. The statement on $\theta_0$-eigenspaces follows by comparing the explicit description from \ref{unimod} with that from \cite{rapo}.\hfill$\Box$\\\\

Of particular interest for the representation theory of ${\rm GL}_{d+1}(k)$ is the virtual ${\rm GL}_{d+1}(k)$-representation $$\sum_{s=0}^d(-1)^sH^s_{rig}(X_{\overline{k}})_{0}=\sum_{s=0}^d(-1)^sH^s_{rig}(Y^0_{\overline{k}})$$ (it is called {\it unipotent}). In view of \ref{bacauni} we determined its reduction modulo $p$ in \ref{unimod}. By (\ref{shife}) (which yields (\ref{spepo})) the situation for the $\theta_j$-eigenspaces of other singular $\theta_j$ is similar.

 To compute the reduction modulo $p$ for non singular $\theta_j$ the general results from section \ref{relri} apply. Alternatively, one might try to understand the finitely generated free $W$-module$$H_j^d=H_{crys}^d(Y/W,{E}(j))/\mbox{torsion}.$$$H_j^d$ is a ${\rm GL}_{d+1}(k)$-stable $W$-lattice in $H^d_{rig}(X_{\overline{k}})_{j}= H_{crys}^d(Y/W,{E}(j))\otimes_{W}K$.
 
\begin{kor} If $\theta_j$ is non singular then $H_j^d\otimes_{W}{k}$ is a ${\rm GL}_{d+1}(k)$-equivariant subquotient of $H^d(Y,(\Omega_{Y}^{\bullet}\otimes{\mathcal L}_{Y}(D(j)),\nabla_j))$ of ${k}$-dimension $(q-1)(q^2-1)\cdots(q^d-1).$ Viewed as an element of the Grothendieck group of $k[{\rm GL}_{d+1}(k)]$-modules it coincides with $$\sum_{s,m}(-1)^{s+m}H^s(Y,\Omega_{Y}^{m}\otimes{\mathcal L}_{Y}(D(j))).$$
\end{kor}

{\sc Proof:} The $k$-dimension of $H_j^d\otimes_{W}{k}$ is the $K$-dimension of $H_j^d\otimes_{W}{K}$ and this can be read off from Corollary \ref{dimfo}. We have a natural surjection $$H_{crys}^d(Y/W,{E}(j))\otimes_{W}{k}\longrightarrow H_j^d\otimes_{W}{k}.$$But $H_{crys}^d(Y/W,{E}(j))\otimes_{W}{k}$ is a subobject of $H_{crys}^d(Y/{k},{E}(j)\otimes_{W}{k})=H^d(Y,(\Omega_{Y}^{\bullet}\otimes{\mathcal L}_{Y}(D(j)),\nabla_j))$, hence the first statement. The second one follows from Theorem \ref{brauer}.\hfill$\Box$\\

{\it Remark:} If $d=1$ and $j$ is arbitrary (but not divisible by $q+1$), or if $d=2$ and $1\le j\le p-1$, one can show $H^s(Y,(\Omega_{Y}^{\bullet}\otimes{\mathcal L}_{Y}(D(j)),\nabla_j))=0$ for $s\ne d$ and by the usual devissage for crystalline cohomology relative to $W/p^n$ for all $n$ we find $H_j^d=H_{crys}^d(Y/W,{E}(j))$ and $$H_j^d\otimes_{W}{k}=H^d(Y,(\Omega_{Y}^{\bullet}\otimes{\mathcal L}_{Y}(D(j)),\nabla_j)).$$The case $d=2$ and $j=p$ is the first one where $H^{d-1}(Y,(\Omega_{Y}^{\bullet}\otimes{\mathcal L}_{Y}(D(j)),\nabla_j))\ne0$ and hence $H_j^d\ne H_{crys}^d(Y/W,{E}(j))$ and $H_j^d\otimes_{W}{k}\ne H^d(Y,(\Omega_{Y}^{\bullet}\otimes{\mathcal L}_{Y}(D(j)),\nabla_j))$ in that case.


\end{document}